\documentclass{article}

\usepackage{amssymb,amsfonts}
\usepackage{amsfonts,amsmath,amsxtra,amsthm,amssymb,latexsym,mathrsfs,srcltx,dsfont,eufrak,euscript}

\usepackage[cp1251]{inputenc}
\usepackage[english]{babel}

\usepackage{hhline}

\newtheorem{defin}{Definition}
\newtheorem{prop}{Proposition}
\newtheorem{theorem}{Theorem}
\newtheorem{lemma}{Lemma}
\newtheorem{conseq}{Corollary}

\newcommand{\AbsEng}[1]{
    \begin{center}
    \begin{minipage}{0.8\linewidth}
      #1
    \end{minipage}
    \end{center}}

\def\R{\mathbb R}

\def\N{\mathbb N}

\def\E{\mathbb E}
\def\Eg{\mathcal E}
\def\G{\mathbb G}

\def\Lg{\mathcal{L}}
\def\F{\mathbb F}

\def\ltri|{\left|\!\!\!\!\!\!\;\;\;\left|\!\!\!\!\!\!\;\;\;\left|}
\def\rtri|{\right|\!\!\!\!\!\!\;\;\;\right|\!\!\!\!\!\!\;\;\;\right|}

\usepackage[usenames]{color}
\usepackage{colortbl}
\usepackage{xcolor}

\usepackage{graphicx}

\begin{document}
\large

\vspace{1cm}

\begin{center} \bf
 A SUPEREXPONENTIALLY CONVERGENT FUNCTIONAL-DISCRETE METHOD FOR SOLVING THE CAUCHY PROBLEM FOR SYSTEMS OF ORDINARY DIFFERENTIAL EQUATIONS\end{center}
\vspace{0.5 cm}
\begin{center} \bf
Makarov V.L.\footnote{E-mail: makarov@imath.kiev.ua} and Dragunov D.V.\footnote{E-mail: dragunovdenis@gmail.com}
\end{center}

\begin{center}\small\it
Department of Numerical Mathematics, \\
Institute of  Mathematics of NAS of Ukraine, \\
3 Tereshchenkivs'ka Str., Kyiv-4, 01601, Ukraine
\end{center}

\vspace{0.25 cm}
\vspace{0.25 cm}

\AbsEng{\small

 In the paper a new numerical-analytical method for solving the Cauchy problem for systems of ordinary differential equations of special form is presented. The method is based on the idea of the FD-method for solving the operator equations of general form, which was proposed by V.L. Makarov. The sufficient conditions for the method converges with a superexponential convergence rate were obtained. We have generalized the known statement about the local properties of Adomian polynomials for scalar functions on the operator case. Using the numerical examples we make the comparison between the proposed method and the Adomian Decomposition Method. }

\vspace{0.5 cm}

{\bf Introduction.}
A great number of papers published during last two decades are devoted to the Adomian Decomposition Method (ADM)  and its applications (see \cite{Adomian_1994} -- \cite{Ercan_Mustafa_DAEs_2006} and the references therein). This method was proposed in the middle of 80-th by American physicist George Adomian. It has aroused great interest among mathematicians all over the world. The very foundations of this interest lie in the fact that this method belongs to the class of analytical methods: it allows one to compute approximate solution of nonlinear operator equation by solving a sequence of linear equations. Using the ADM we looking for the approximate solution in the form of partial sum of the series which converges to the exact solution. Due to its analytical property the ADM allows one to obtain the approximate solution even in the case when initial operator equation depends on the parameters. In other words, the ADM allows us to approximate the multi-parametric family of solutions in the functional form.

The general idea of the functional-discrete method (or simply FD-method) is quite similar to the ADM's one. The FD-method was for the first time proposed in \cite{Makarov_FD_1991} where it was applied  to the Sturm–Liouville eigenvalue problem. In \cite{Makarov_FD_1991} it was showed that the convergence rate of the method increases as far as increases the  order number of the eigenvalue which we are going to approximate. So, the FD-method demonstrated amazing results that are unachievable for the finite elements method, mesh method and other purely discrete methods. Leter, in \cite{G_K_M_R_2007} -- \cite{Makarov_DAN_2008_FD_Vl_Zn}, the FD-method was applied to the nonlinear eigenvalue problems and  the sufficient conditions that guarantee its superexponential convergence rate were found. Applying the FD-method to the boundary value problem, the authors of \cite{Gav_Laz_Mak_Sytn_2008} proved, that in this case the superexponential convergence rate of the method can be reached (under certain sufficient conditions) as well. In addition to this, in \cite{Gav_Laz_Mak_Sytn_2008} the general scheme of the FD-method for some class of operator equations in an abstract normed spaces was presented. In \cite{Drag_4} the FD-method was applied to the Cauchy problem for ordinary differential equation of the first order and the sufficient conditions which guarantee superexponential convergence rate of the method were found.

Both the ADM and the FD-method are originated from the homotopy method. Both this methods approximate the exact solution of the nonlinear problem using the partial sums, which terms are to be found as the solutions of some sequence of linear problems\footnote{Generally speaking, this is not true for the FD-method, because the base problem is not linear one (see later).}. But the essential difference between those methods is expressed by the fact that the FD-method has a built-in parameter --- the parameter of discretization. By varying this parameter we cane achieve the convergence of the FD-method even in the cases when the ADM is turn out to be divergent. The last fact will be demonstrated below by several numerical examples.

In the present paper we offer the algorithm of the FD-method for solving the Cauchy problem for systems of ordinary differential equations (SODE), which are based on the general scheme of the FD-method for operator equations \cite{Gav_Laz_Mak_Sytn_2008}. It is the further development and generalization of the theoretical results obtained in \cite{Mak_Drag_1}. The main result of the paper is presented in theorem \ref{T_first_FD}, which contains the sufficient conditions for the FD-method applied to SODE converges with superexponential rate.

Let us consider the following Cauchy problem for SODE of the first order
\begin{equation}\label{Chapt_3_1_eq_1}
    \frac{d}{d t}\overrightarrow{u}\left(t\right) - \mathbf{N}\left(t, \overrightarrow{u}\left(t\right)\right)\overrightarrow{u}\left(t\right)=\overrightarrow{\phi}\left(t\right),\; \overrightarrow{u}\left(t_{0}\right)=\overrightarrow{u}_{0}\in V_{m}\left(\R\right),\quad t\in\left[t_{0}, +\infty\right),
\end{equation}
 where $\mathbf{N}\left(t, \overrightarrow{u}\right)\in M_{m}\left(C_{t,\overrightarrow{u}}^{0,\infty}\left(\left[t_{0}, +\infty\right)\times \R^{m}\right)\right),$ $\overrightarrow{u}\left(t\right)=\left[u_{1}\left(t\right),\ldots, u_{m}\left(t\right)\right]^{T},$ $\overrightarrow{\phi}\left(t\right)\in V_{m}\left(C\left(\left[t_{0}, +\infty\right)\right)\right).$ Here and below we use the notation $M_{m}\left(\F\right)$ to describe a linear space of square matrices of the $m$-th order with elements belonging to a linear space $\F.$  Similarly to this, by  $V_{m}\left(\F\right)$ we denote a linear space of column vectors of the $m$-th order with elements from $\F.$ We also assume that problem (\ref{Chapt_3_1_eq_1}) possess a unique solution.

Using the FD-method for solving the Cauchy problem (\ref{Chapt_3_1_eq_1}), we looking for the approximation of the exact solution $\overrightarrow{u}\left(t\right)$ in the form of the partial sum (see \cite{Gav_Laz_Mak_Sytn_2008})
\begin{equation}\label{Chapt_3_1_eq_2}
    \overset{p}{\overrightarrow{u}}\left(t\right)=\sum\limits_{i=0}^{p}\overrightarrow{u}^{(i)}\left(t\right),\; p\in\N\bigcup\left\{0\right\}
\end{equation}
of the convergent series
\begin{equation}\label{Chapt_3_1_eq_3}
    \overrightarrow{u}\left(t\right)=\sum\limits_{i=0}^{\infty}\overrightarrow{u}^{(i)}\left(t\right).
\end{equation}
In this case the nonnegative  integer number $p$ is called the rank of the FD-method. To emphasize this detail we usually say:``the FD-method of the $p$-th rank''.

To apply the FD-method to the Cauchy problem (\ref{Chapt_3_1_eq_1}) we need to introduce a grid on the interval $\left[t_{0}, +\infty\right)$ :
\begin{equation}\label{Chapt_3_1_eq_4}
\begin{array}{c}
  \widehat{\omega}=\bigg\{t_{0}<t_{1}<t_{2}<\ldots,\;\lim\limits_{n\rightarrow +\infty}t_{n}=+\infty\bigg\},\quad h=\sup\limits_{i\in \N}\left\{h_{i}=t_{i}-t_{i-1}\right\}. \\[1.2em]
\end{array}
\end{equation}
After that we can consider the following recursive system of linear Cauchy problems with respect to the unknown terms $\overrightarrow{u}^{(j)}\left(t\right)=\left[u_{1}^{(j)}\left(t\right),\ldots, u_{m}^{(j)}\left(t\right)\right]^{T},\; j\in\N\bigcup \left\{0\right\}$ of the series (\ref{Chapt_3_1_eq_3}):
\begin{equation}\label{Chapt_3_1_eq_5}
\begin{array}{c}
  \cfrac{d}{d t}\overrightarrow{u}^{(0)}\left(t\right)-\mathbf{N}\left(t, \overrightarrow{u}^{(0)}\left(t_{i-1}\right)\right)\overrightarrow{u}^{(0)}\left(t\right)=\overrightarrow{\phi}\left(t\right),\quad t\in\left[t_{i-1}, t_{i}\right), \forall i\in \N, \\ [1.2em]
  \overrightarrow{u}^{(0)}\left(t_{0}\right)=\overrightarrow{u}_{0},\quad \left[\overrightarrow{u}^{(0)}\left(t\right)\right]_{t=t_{i}}=\overrightarrow{u}^{(0)}\left(t_{i}+0\right)-\overrightarrow{u}^{(0)}\left(t_{i}-0\right)=\overrightarrow{0},\; \forall i\in \N.
\end{array}
\end{equation}
\begin{equation}\label{Chapt_3_1_eq_6}
    \begin{array}{c}
      \cfrac{d}{dx}\overrightarrow{u}^{(j)}\left(t\right)-\mathbf{N}\left(t, \overrightarrow{u}^{(0)}\left(t_{i-1}\right)\right)\overrightarrow{u}^{(j)}\left(t\right)= \\[1.2em]
      =\left[\sum\limits_{p=1}^{m}\left(\cfrac{\partial }{\partial
                            u_{p}}\mathbf{N}\left(t, \overrightarrow{u}^{(0)}\left(t_{i-1}\right)\right)\right)\;u_{p}^{(j)}\left(t_{i-1}\right)\right]\overrightarrow{u}^{(0)}\left(t\right)+F^{(j)}\left(t\right),\;t\in\left[t_{i-1}, t_{i}\right),\\[1.2em]
      \overrightarrow{u}^{(j)}\left(t_{0}\right)=\overrightarrow{0},\quad \left[\overrightarrow{u}^{(j)}\left(t\right)\right]_{t=t_{i}}=\overrightarrow{0},\; \forall i,j\in \N,
    \end{array}
\end{equation}
where
\begin{equation}\label{Chapt_3_1_eq_7}
   \begin{array}{c}
     F^{(j)}\left(t\right)=\sum\limits_{p=1}^{j-1}A_{j-p}\left(\mathbf{N}\left(t, \left(\cdot\right)\right); \overrightarrow{u}^{(0)}\left(t_{i-1}\right), \ldots , \overrightarrow{u}^{(j-p)}\left(t_{i-1}\right)\right)\overrightarrow{u}^{(p)}\left(t\right)+ \\[1.2em]
     +\sum\limits_{p=0}^{j-1}\left[A_{j-1-p}\left(\mathbf{N}\left(t, \left(\cdot\right)\right); \overrightarrow{u}^{(0)}\left(t\right), \ldots, \overrightarrow{u}^{(j-1-p)}\left(t\right) \right)-\right. \\[1.2em]
     \left.-A_{j-1-p}\left(\mathbf{N}\left(t, \left(\cdot\right)\right); \overrightarrow{u}^{(0)}\left(t_{i-1}\right), \ldots, \overrightarrow{u}^{(j-1-p)}\left(t_{i-1}\right) \right)\right]\overrightarrow{u}^{(p)}\left(t\right)+ \\[1.2em]
     +A_{j}\left(\mathbf{N}\left(t, \left(\cdot\right)\right); \overrightarrow{u}^{(0)}\left(t_{i-1}\right), \ldots, \overrightarrow{u}^{(j-1)}\left(t_{i-1}\right), \overrightarrow{0}\right)\overrightarrow{u}^{(0)}\left(t\right),\\[1.2em]
     t\in \left[t_{i-1}, t_{i}\right),\quad \forall i,j\in \N.
   \end{array}
\end{equation}
Here $A_{k}\left(\mathbf{N}\left(t, \left(\cdot\right)\right); \overrightarrow{v}_{0},\overrightarrow{v}_{1},\ldots, \overrightarrow{v}_{k},\right)=\left.\frac{1}{k!}\frac{d^{k}}{d \tau^{k}}\mathbf{N}\left(t, \sum\limits_{i=0}^{\infty}\tau^{i}\overrightarrow{v}_{i}\right)\right|_{\tau=0}$ denotes the Adomian's polynomial of the $k\in \N\bigcup\left\{0\right\}$ order for the operator $\mathbf{N}\left(t, \left(\cdot\right)\right)$ ( see \cite{Seng_Abbaoui_Cherruault}, \cite{Seng_Abbaoui_Cherruault_1}).

 Problem (\ref{Chapt_3_1_eq_5}) is called \textit{the base problem}. It can be considered as an example of the Cauchy problem with \textit{piecewise constant argument}   (see \cite{Akhmet_PC_argument}).
If we denote by $U_{i}\left(t\right)=U_{i}\left(t,\;\mathbf{N}\left(t, \overrightarrow{u}^{(0)}\left(t_{i-1}\right)\right)\right)\in M_{m}\left(C^{1}\left[t_{i-1}, t_{i}\right]\right),$ $i\in \N$ the solution of the following Cauchy problem
\begin{equation}\label{Chapt_3_1_eq_8}
\begin{array}{c}
  \frac{d}{d t}U_{i}\left(t\right)-\mathbf{N}\left(t, \overrightarrow{u}^{(0)}\left(t_{i-1}\right)\right) U_{i}\left(t\right)=0, \quad U_{i}\left(t_{i-1}\right)=E,\quad t\in\left[t_{i-1}, t_{i}\right],
\end{array}
\end{equation}
the solution of problem (\ref{Chapt_3_1_eq_5}) can be expressed in the form
\begin{equation}\label{Chapt_3_1_eq_10}
    \overrightarrow{u}^{(0)}\left(t\right)=U_{i}\left(t\right)\overrightarrow{u}^{(0)}\left(t_{i-1}\right)+\int\limits_{t_{i-1}}^{t}K_{i}\left(t, \xi\right)\overrightarrow{\phi}\left(\xi\right)d \xi,\quad t\in\left[t_{i-1}, t_{i}\right],\; \forall i\in \N,
\end{equation}
where
\begin{equation}\label{Chapt_3_1_eq_11}
    K_{i}\left(t, \xi\right)=U_{i}\left(t\right)U_{i}^{-1}\left(\xi\right),\quad t\in\left[t_{i-1}, t_{i}\right],\; \forall i\in \N,
\end{equation}
is the Cauchy matrix (see, for example, \cite[p. 412]{Gantmaher}).
Similarly we can express the solutions of problems (\ref{Chapt_3_1_eq_6}), (\ref{Chapt_3_1_eq_7}) $\forall j\in\N$:
\begin{equation}\label{Chapt_3_1_eq_12}
\begin{array}{c}
  \overrightarrow{u}^{(j)}\left(t\right)=U_{i}\left(t\right)\overrightarrow{u}^{(j)}\left(t_{i-1}\right)+\int\limits_{t_{i-1}}^{t}K\left(t, \xi\right)\Upsilon_{i}\left(\xi\right)\overrightarrow{u}^{(j)}\left(t_{i-1}\right) d\xi+ \\[1.2em]
  +\int\limits_{t_{i-1}}^{t}K\left(t, \xi\right) F^{(j)}\left(\xi\right) d\xi,\quad t\in\left[t_{i-1}, t_{i}\right], \forall i\in \N,
\end{array}
\end{equation}
 where $\Upsilon_{i}\left(t\right)\in M_{m}\left(C^{1}\left[t_{i-1}, t_{i}\right]\right)$ denotes the following matrix
\begin{equation}\label{Chapt_3_1_eq_13}
    \Upsilon_{i}\left(t\right)=\left[\frac{\partial }{\partial
                            u_{1}}\mathbf{N}\left(t, \overrightarrow{u}^{(0)}\left(t_{i-1}\right)\right)\overrightarrow{u}^{(0)}\left(t\right),\ldots, \frac{\partial }{\partial
                            u_{m}}\mathbf{N}\left(t, \overrightarrow{u}^{(0)}\left(t_{i-1}\right)\right)\overrightarrow{u}^{(0)}\left(t\right)\right],
\end{equation}
$\forall t\in\left[t_{i-1}, t_{i}\right].$ From expressions (\ref{Chapt_3_1_eq_10}), (\ref{Chapt_3_1_eq_12}) and (\ref{Chapt_3_1_eq_7}) we can easily obtain the following statement.
\begin{prop}
  The column vectors $\overrightarrow{u}^{(0)}\left(t\right), \overrightarrow{u}^{(1)}\left(t\right), \ldots,$ which are the solutions of the recursive  linear Cauchy problems \eqref{Chapt_3_1_eq_5}--\eqref{Chapt_3_1_eq_7}, always exist and are unique on $\left[t_{0}, +\infty\right].$
\end{prop}

\begin{defin}\label{O_FD_zbizhn}
   We will say that the FD-method for the Cauchy problem \eqref{Chapt_3_1_eq_1} converges {\normalfont (}to the exact solution of problem \eqref{Chapt_3_1_eq_1}{\normalfont )} on $\left[t_{0}, t_{0}+H\right),$ \mbox{$t_{0}<H\leq +\infty,$} if there exists a positive constant $\overline{h}\in \R$ such that for every grid $\widehat{\omega}$ \eqref{Chapt_3_1_eq_4} satisfying  $h\leq \overline{h},$  series \eqref{Chapt_3_1_eq_3}, with terms obtained as the solutions of problems \eqref{Chapt_3_1_eq_5} -- \eqref{Chapt_3_1_eq_7}, converges absolutely and uniformly  {\normalfont (}to the exact solution of problem  \eqref{Chapt_3_1_eq_1}{\normalfont )} on $\left[t_{0}, t_{0}+H\right).$
\end{defin}

{\bf The local properties of the Adomian's polynomials.}
  Let $\E_{1}$ and $\E_{2}$ be some normed linear spaces over $\R$ with norms $\left\|\cdot\right\|^{(1)}$ and $\left\|\cdot\right\|^{(2)}$ respectively. A set of linear operators which map $\E_{1}$ to (on) $\E_{2}$ we will denote by $\Lg\left(\E_{1}, \E_{2}\right).$
  Lets consider an operator $\mathbf{N}\left(\mathbf{u}\right)$ (nonlinear, in general), which maps a given open subset $\G$ of the space $\E_{1}$ into $\E_{2},$ i.e. $\mathbf{N}\colon \G\rightarrow \E_{2}.$  Let us recall some definitions of important concepts, which we will need later   ( see \cite[p. 28--30]{Mak_Hlob_Yan_IO}).
  \begin{defin}
     An operator $\mathbf{N}^{(k)}\left(\mathbf{u}_{1}, \mathbf{u}_{2},\ldots, \mathbf{u}_{k}\right)\colon \E_{1}^{k}\rightarrow \E_{2}$ is said to be a \textbf{$k$-linear operator}, if it is linear with respect to each argument $\mathbf{u}_{i}\in\E_{1},$ $i\in\overline{1, k}.$

     An operator  $\mathbf{N}^{(k)}\left(\mathbf{u}_{1}, \mathbf{u}_{2},\ldots, \mathbf{u}_{k}\right)$ is called to be \textbf{symmetric}, if its value is independent on the order of its arguments.
  \end{defin}
  We will denote the set of $k$-linear operators $\mathbf{N}^{(k)}\left(\mathbf{u}_{1}, \mathbf{u}_{2},\ldots, \mathbf{u}_{k}\right)\colon \E_{1}^{k}\rightarrow \E_{2}$ by $\Lg_{k}\left(\E_{1}^{k}, \E_{2}\right).$
  \begin{defin}\label{O_normi_klin_operatora}
     A given $k$-linear operator $\mathbf{N}^{(k)}\left(\mathbf{u}_{1}, \mathbf{u}_{2},\ldots, \mathbf{u}_{k}\right)$ is said to be \textbf{bounded}, if there exists a constant $M$ such that for every $\mathbf{u}_{i}\in \E_{1}$ the following inequality holds
     \begin{equation}\label{operat_lema_15}
        \left\|\mathbf{N}^{(k)}\left(\mathbf{u}_{1}, \mathbf{u}_{2},\ldots, \mathbf{u}_{k}\right)\right\|^{(2)}\leq M\left\|\mathbf{u}_{1}\right\|^{(1)}\left\|\mathbf{u}_{2}\right\|^{(1)}\ldots \left\|\mathbf{u}_{k}\right\|^{(1)}.
     \end{equation}
     For the smallest constant $M$ which satisfies inequality \eqref{operat_lema_15} we will use notation $\left\|\mathbf{N}^{(k)}\right\|.$ It is called the \textbf{norm} of the $k$-linear operator $\mathbf{N}^{(k)}.$ In other words,
     \begin{equation}\label{operat_lema_16}
        \left\|\mathbf{N}^{(k)}\right\|=\sup\limits_{\left\|\mathbf{u}_{1}\right\|^{(1)}\leq 1,\ldots, \left\|\mathbf{u}_{k}\right\|^{(1)}\leq 1}\left\|\mathbf{N}^{(k)}\left(\mathbf{u}_{1}, \mathbf{u}_{2},\ldots, \mathbf{u}_{k}\right)\right\|^{(2)}.
     \end{equation}
  \end{defin}
 \begin{defin}\label{O_pohidnoi_Freshet}
 An operator $\mathbf{N}\left(\mathbf{u}\right)$ is said to be \textbf{differentiable by the Frechet} at the point $\mathbf{u}_{0}\in\G$ if there exists a bounded linear operator $\mathbf{N}^{(1)}\left(\mathbf{u}_{0}\right)\in \Lg\left(\E_{1}, \E_{2}\right)$ such that
 \begin{equation}\label{operat_lema_17}
    \mathbf{N}\left(\mathbf{u}_{0}+\mathbf{h}\right)-\mathbf{N}\left(\mathbf{u}_{0}\right)=\mathbf{N}^{(1)}\left(\mathbf{u}_{0}\right)\mathbf{h}+\mathbf{\omega}\left(\mathbf{u}_{0}, \mathbf{h}\right),
 \end{equation}
  where $\left\|\mathbf{\omega}\left(\mathbf{u}_{0}, \mathbf{h}\right)\right\|^{(2)}=o\left(\left\|\mathbf{h}\right\|^{(1)}\right),$ $\left\|\mathbf{h}\right\|^{(1)}\rightarrow 0,$ $\forall \mathbf{h}\colon \mathbf{u}_{0}+\mathbf{h}\in \G.$ The linear operator $\mathbf{N}^{(1)}\left(\mathbf{u}_{0}\right),$ when it exists, is called the 1-st \textbf{Frechet derivative} of the operator $\mathbf{N}\left(\mathbf{u}\right)$ at the point $\mathbf{u}_{0}.$
 \end{defin}

Using the induction we can define the $k$-th Frechet derivative of the operator $\mathbf{N}\left(\mathbf{u}\right)$ at the point $\mathbf{u}_{0},$ $\forall k\in\N$ \cite[p. 262]{shvartz_analiz}. Namely, if the 1-st Frechet derivative of the operator $\mathbf{N}\left(\mathbf{u}\right)$ exists at each point of the certain open neighbourhood  $\G_{1}\subseteq \G$ of the point $\mathbf{u}_{0},$ then we can consider the mapping (or the operator)
\begin{equation}\label{operat_lema_20}
 \mathbf{N}^{(1)}\left(\mathbf{u}\right) \colon \G_{1} \rightarrow \Lg\left(\E_{1}, \E_{2}\right).
\end{equation}
 If there exists the 1-st Frechet derivative of operator (\ref{operat_lema_20}) at the point $\mathbf{u}_{0}$, then we will denote it by $\mathbf{N}^{(2)}\left(\mathbf{u}_{0}\right)$ and call it the second Frechet derivative of the operator $\mathbf{N}\left(\mathbf{u}\right)$ at the point $\mathbf{u}_{0}.$ It is obvious that $\mathbf{N}^{(2)}\left(\mathbf{u}_{0}\right)\in\Lg_{2}\left(\E_{1}^{2}, \E_{2}\right).$ In general case, the $k$-th Frechet derivative of the operator $\mathbf{N}\left(\mathbf{u}\right)$ at the point $\mathbf{u}_{0}$ belongs to the space $\Lg_{k}\left(\E^{k}_{1}, \E_{2}\right).$ It is denoted by $\mathbf{N}^{(k)}\left(\mathbf{u}_{0}\right)$ and is defined as the 1-st Frechet derivative of the operator $\mathbf{N}^{(k-1)}\left(\mathbf{u}\right) \colon \G_{k-1} \rightarrow \Lg_{k-1}\left(\E^{k-1}_{1}, \E_{2}\right)$ at the point $\mathbf{u}_{0}\in \G_{k-1}.$
\begin{defin}
  An operator $\mathbf{N}\left(u\right)$ is said to be $k$ {\normalfont (}$k\in\N${\normalfont )} times differentiable by the Frechet (or in the Frechet's sense) on the certain open set $\G_{k}\subseteq \G,$ if  there exists the $k$-th Frechet derivative of the operator $\mathbf{N}\left(u\right)$ at each point of $\G_{k}.$
\end{defin}

Below we give the definition of the Adomian polynomial of the $n$-th order for the operator $\mathbf{N}\left(\mathbf{u}\right)$ with respect to the ``variables'' $\mathbf{u}_{i}\in \E,$ $i\in \overline{1,n}.$ This definition is more detailed then the definition from \cite[p. 59]{Seng_Abbaoui_Cherruault}.
\begin{defin}\label{O_oper_poly_adom}
      Let the operator $\mathbf{N}\left(\mathbf{u}\right)\colon \G\rightarrow \E_{2}$ is differentiable by the Frechet at the point $\mathbf{u}_{0}\in \G$ up to the $n$-th order inclusively, $n\in\N.$ Additionally, let $\Phi\left(\tau\right)\overset{{\normalfont \textmd{def}}}{=}\Phi\left(\tau, \mathbf{u}_{1},\ldots,\mathbf{u}_{n}\right)$ be an operator which maps  $\R\times \E_{1}^{n}$  {\normalfont (}$\tau \in \mathbb{R},$ $\mathbf{u}_{i}\in \E_{1}$, $i\in\overline{1, n}${\normalfont )} to $\E_{1}$ and is defined by virtue of the formula $\Phi\left(\tau\right)=\sum\limits_{i=0}^{n}\tau^{i} \mathbf{u}_{i}.$ Then the Adomian polynomial of the $n$-th order for the operator $\mathbf{N}\left(\mathbf{u}\right)$ is defined to be an operator $A_{n}\left(\mathbf{N}\left(\cdot\right); \left[\mathbf{u}_{i}\right]_{i=0}^{n}\right)\overset{{\normalfont \textmd{def}}}{=}A_{n}\left(\mathbf{N}\left(\cdot\right); \mathbf{u}_{0}, \mathbf{u}_{1}, \ldots, \mathbf{u}_{n}\right),$ which maps $\E_{1}^{n}$ to $E_{2}$ and can be explicitly obtained by the formula
   \begin{equation}\label{operator_lema_2}
    A_{n}\left(\mathbf{N}\left(\cdot\right); \left[\mathbf{u}_{i}\right]_{i=0}^{n}\right)=\left.\cfrac{1}{n!}\cfrac{d^{n}}{d \tau^{n}}\mathbf{N}\left(\Phi\left(\tau\right)\right)\right|_{\tau=0}.
   \end{equation}
    If $n=0$ then we assume $$A_{0}\left(\mathbf{N}\left(\cdot\right); \mathbf{u}_{0}\right)\equiv \mathbf{N}\left(\mathbf{u}_{0}\right).$$
\end{defin}

 For convenience we will use the following notation for the  $k$-linear operator $\mathbf{N}^{(k)}\left(\mathbf{u}_{0}\right)\in \Lg_{k}\left(\E_{1}^{k}, \E_{2}\right)$
\begin{equation}\label{operat_lema_7}
\begin{array}{c}
  \mathbf{N}^{(k)}\left(\mathbf{u}_{0}\right)\bigg(\overbrace{\mathbf{u}_{1},\ldots, \mathbf{u}_{1}}\limits^{k_{1}},\overbrace{\mathbf{u}_{2},\ldots, \mathbf{u}_{2}}\limits^{k_{2}},\ldots,\overbrace{\mathbf{u}_{l},\ldots, \mathbf{u}_{l}}\limits^{k_{l}} \bigg)\overset{\textmd{def}}{=} \\[1.2em]
  \overset{\textmd{def}}{=}\mathbf{N}^{(k)}\left(\mathbf{u}_{0}\right)\left(\overset{k_{1}}{\left[\mathbf{u}_{1}\right]},\overset{k_{2}}{\left[\mathbf{u}_{2}\right]},\ldots,\overset{k_{l}}{\left[\mathbf{u}_{l}\right]} \right),\\[1.2em]
  \forall \mathbf{u}_{i}\in \E_{1},\quad k_{i}\in\N\bigcup \left\{0\right\},\quad i\in \overline{1, l},\quad \sum\limits_{j=1}^{l}k_{j}=k.
\end{array}
\end{equation}

\begin{theorem}{\normalfont (see \cite[p. 60]{Seng_Abbaoui_Cherruault})}\label{T_pro_O_M_A}
Let an operator $\mathbf{N}\left(u\right)$ is differentiable by the Frechet at the point $\mathbf{u}_{0}\in\G$ up to the $n$-th order inclusively, $n\in \N,$ then $\forall \mathbf{u}_{i}\in \E_{1},$ $\forall i\in \overline{1,n}$
\begin{equation}\label{operat_lema_8}
  A_{n}\left(\mathbf{N}\left(\cdot\right); \left[\mathbf{u}_{i}\right]_{i=0}^{n}\right)=\sum\limits_{\substack{\sum\limits_{k=1}^{n}k p_{k}=n \\ p_{k}\in \N\bigcup \left\{0\right\}}}\cfrac{1}{p_{1}!\ldots p_{n}!}\mathbf{N}^{\left(p_{1}+\ldots+p_{n}\right)}\left(\mathbf{u}_{0}\right)\left(\overset{p_{1}}{\left[\mathbf{u}_{1}\right]},\ldots, \overset{p_{n}}{\left[\mathbf{u}_{n}\right]}\right),
\end{equation}
\begin{equation}\label{operat_lema_8'}
\begin{array}{c}
  A_{n}\left(\mathbf{N}\left(\cdot\right); \left[\mathbf{u}_{i}\right]_{i=0}^{n}\right)=\sum\limits_{\substack{ \alpha_{1}+\ldots+\alpha_{n}=n \\ \alpha_{1}\geq\ldots\geq\alpha_{n+1}=0 \\ \alpha_{i}\in \N\bigcup \left\{0\right\}}}\cfrac{1}{\prod\limits_{i=1}^{n}\left(\alpha_{i}-\alpha_{i+1}\right)!}\mathbf{N}^{\left(\alpha_{1}\right)}\left(\mathbf{u}_{0}\right)\left(\overset{\alpha_{1}-\alpha_{2}}{\left[\mathbf{u}_{1}\right]},\ldots, \overset{\alpha_{n}-\alpha_{n+1}}{\left[\mathbf{u}_{n}\right]}\right).
\end{array}
\end{equation}

\end{theorem}
Lets consider the subspace $\Eg_{1}$ of the normed space $\E_{1},$ which is also a normed space with the norm $\left\|\cdot\right\|^{(1)}$.
Let $\left\|\cdot\right\|^{(1)}_{1}$ be a scalar function defined on  $\Eg_{1},$ which satisfying the conditions of seminorm (see, for example, \cite[p. 564]{shvartz_analiz}), then a scalar function $\ltri|\cdot\rtri|,$ defined by the equality
\begin{equation}\label{operat_lema_25}
    \ltri|\mathbf{u}\rtri|=\max\left\{\left\|\mathbf{u}\right\|^{(1)}, \left\|\mathbf{u}\right\|^{(1)}_{1}\right\},\quad \forall \mathbf{u}\in \Eg_{1},
\end{equation}
is a norm (see, for example, \cite[p. 41--42]{shvartz_analiz}) in the space $\Eg_{1}.$

The following lemma is a generalization of lemma 2.1 from \cite{Gav_Laz_Mak_Sytn_2008}.
\begin{lemma}\label{L_pro_Adom_1_Operat}
Let the following conditions are satisfied
\begin{enumerate}
\item An operator $\mathbf{N}\left(u\right)$ is $n$ times {\normalfont (}$n\in\N${\normalfont )} differentiable by the Frechet on the open convex subspace $\G_{n}$ of the space $\E_{1},$ $\G_{n}\bigcap \Eg_{1}\neq \varnothing$ and there exists a scalar function $\widetilde{N}\left(u\right)\in C^{n}\left(\mathbb{R}\right)$ which satisfies the following inequalities
\begin{equation}\label{operat_lema_21}
    \left\|\mathbf{N}^{(k)}\left(\mathbf{u}\right)\right\|\leq \left.\cfrac{d^{k}}{d u^{k}}\widetilde{N}\left(u\right)\right|_{u=\left\|\mathbf{u}\right\|^{(1)}},\quad \forall k\in\overline{0, n},\quad \forall \mathbf{u}\in\G_{n}.
\end{equation}
\item For some positive constant $h\in \mathbb{R}$ and operator $\mathbf{P}_{h}\colon \Eg_{1}\rightarrow \E_{1}$ the following conditions hold true
\begin{equation}\label{operat_lema_22}
    \left\|\mathbf{u}-\mathbf{P}_{h}\left(\mathbf{u}\right)\right\|^{(1)}\leq h\left\|\mathbf{u}\right\|^{(1)}_{1}, \quad \forall \mathbf{u}\in\Eg_{1},
\end{equation}
\begin{equation}\label{operat_lema_22'}
    \mathbf{P}_{h}\left(\mathbf{u}\right)\in \G_{n},\quad \forall \mathbf{u}\in \G_{n}\bigcap \Eg_{1},
\end{equation}
\begin{equation}\label{operat_lema_22''}
\begin{array}{c}
  \left\|\mathbf{P}_{h}\left(\mathbf{u}\right)-\theta_{1}\left(\mathbf{P}_{h}\left(\mathbf{u}\right)-\mathbf{u}\right)\right\|^{(1)}\geq \left\|\mathbf{P}_{h}\left(\mathbf{u}\right)-\theta_{2}\left(\mathbf{P}_{h}\left(\mathbf{u}\right)-\mathbf{u}\right)\right\|^{(1)}, \\[1.2em]
  \quad \forall \theta_{1}, \theta_{2}\in \mathbb{R} \colon \theta_{1}\geq \theta_{2}\geq 0,\quad\forall \mathbf{u}\in \Eg_{1}.
\end{array}
\end{equation}
\end{enumerate}
Then
\begin{equation}\label{operat_lema_23}
\begin{array}{c}
  \left\|A_{k}\left(\mathbf{N}\left(\cdot\right); \left[\mathbf{u}_{i}\right]_{i=0}^{k}\right)\right\|^{(2)}\leq A_{k}\left(\widetilde{N}\left(\cdot\right); \left[\left\|\mathbf{u}_{i}\right\|^{(1)}\right]_{i=0}^{k}\right),  \\[1.2em]
   \forall \mathbf{u}_{0}\in\G_{n},\quad \forall \mathbf{u}_{i}\in \E_{1},\quad i\in \overline{1,k},\quad \forall k\in\overline{0,n-1,}
\end{array}
\end{equation}
\begin{equation}\label{operat_lema_24}
\begin{array}{c}
  \left\|A_{k}\left(\mathbf{N}\left(\cdot\right); \left[\mathbf{u}_{i}\right]_{i=0}^{k}\right)-A_{k}\left(\mathbf{N}\left(\cdot\right); \left[\mathbf{P}_{h}\left(\mathbf{u}_{i}\right)\right]_{i=0}^{k}\right)\right\|^{(2)}\leq \\ [1.2em]
  \leq h A_{k}\left(\widetilde{N}^{(1)}\left(\cdot\right)\times\left(\cdot\right); \left[\ltri|\mathbf{u}_{i}\rtri|\right]_{i=0}^{k}\right),\\[1.2em]
  \forall \mathbf{u}_{0}\in \G_{n}\bigcup \Eg_{1},\quad \forall \mathbf{u}_{i}\in \Eg_{1},\quad i\in \overline{1,k},\quad \forall k\in\overline{0,n-2}.
\end{array}
\end{equation}
\end{lemma}
{\bf Proof.} Let the conditions of the lemma are fulfilled. At a first step we are going to prove inequalities (\ref{operat_lema_23}). Lets fix an arbitrary $k\in \overline{0, n-1}.$ Using theorem \ref{T_pro_O_M_A}, specifically, representation (\ref{operat_lema_8'}), as well as inequalities (\ref{operat_lema_21}), we would obtain
\begin{equation}\label{operat_lema_26}
\begin{array}{c}
  \left\|A_{k}\left(\mathbf{N}\left(\cdot\right); \left[\mathbf{u}_{i}\right]_{i=0}^{k}\right)\right\|^{(2)}\leq\sum\limits_{\substack{ \alpha_{1}+\ldots+\alpha_{k}=k \\ \alpha_{1}\geq\ldots\geq\alpha_{k+1}=0 \\ \alpha_{i}\in \N\bigcup \left\{0\right\}}}\cfrac{\left\|\mathbf{N}^{\left(\alpha_{1}\right)}\left(\mathbf{u}_{0}\right)\right\|\prod\limits_{i=1}^{k}\left(\left\|\mathbf{u}_{i}\right\|^{(1)}\right)^{\alpha_{i}-\alpha_{i+1}}}{\prod\limits_{i=1}^{k}\left(\alpha_{i}-\alpha_{i+1}\right)!}\leq\\[1.2em]
  \leq \sum\limits_{\substack{ \alpha_{1}+\ldots+\alpha_{k}=k \\ \alpha_{1}\geq\ldots\geq\alpha_{k+1}=0 \\ \alpha_{i}\in \N\bigcup \left\{0\right\}}}\cfrac{\widetilde{N}^{\left(\alpha_{1}\right)}\left(\left\|\mathbf{u}_{0}\right\|^{(1)}\right)\prod\limits_{i=1}^{k}\left(\left\|\mathbf{u}_{i}\right\|^{(1)}\right)^{\alpha_{i}-\alpha_{i+1}}}{\prod\limits_{i=1}^{k}\left(\alpha_{i}-\alpha_{i+1}\right)!}=
\end{array}
\end{equation}
$$=A_{k}\left(\widetilde{N}\left(\cdot\right); \left[\left\|\mathbf{u}_{i}\right\|^{(1)}\right]_{i=0}^{k}\right),\quad \forall \mathbf{u}_{0}\in \G_{n}, \quad \forall \mathbf{u}_{i}\in \E_{1},\;i\in \overline{1,k}.$$

Let us prove inequalities (\ref{operat_lema_24}). For this reason we fix an arbitrary $k\in \overline{0, n-2}$ again.  Using representation (\ref{operat_lema_8'}) and the generalization of the Mean value theorem on the operator case (see, for example, \cite[p. 248]{shvartz_analiz}), as well as the properties of the operator's $\mathbf{N}\left(u\right)$ derivatives, we will obtain:
\begin{equation}\label{operat_lema_27}
    \begin{array}{c}
      A_{k}\left(\mathbf{N}\left(\cdot\right); \left[\mathbf{u}_{i}\right]_{i=0}^{k}\right)-A_{k}\left(\mathbf{N}\left(\cdot\right); \left[\mathbf{P}_{h}\left(\mathbf{u}_{i}\right)\right]_{i=0}^{k}\right)= \\[1.2em]
      =\sum\limits_{\substack{ \alpha_{1}+\ldots+\alpha_{k}=k \\ \alpha_{1}\geq\ldots\geq\alpha_{k+1}=0 \\ \alpha_{i}\in \N\bigcup \left\{0\right\}}}\left[\cfrac{\mathbf{N}^{\left(\alpha_{1}\right)}\left(\mathbf{u}_{0}\right)\left(\overset{\alpha_{1}-\alpha_{2}}{\left[\mathbf{u}_{1}\right]},\ldots, \overset{\alpha_{k}-\alpha_{k+1}}{\left[\mathbf{u}_{k}\right]}\right)}{\prod\limits_{i=1}^{k}\left(\alpha_{i}-\alpha_{i+1}\right)!}-\right.\\[1.2em]
      \left.-\cfrac{\mathbf{N}^{\left(\alpha_{1}\right)}\left(\mathbf{P}_{h}\left(\mathbf{u}_{0}\right)\right)\left(\overset{\alpha_{1}-\alpha_{2}}{\left[\mathbf{P}_{h}\left(\mathbf{u}_{1}\right)\right]},\ldots, \overset{\alpha_{k}-\alpha_{k+1}}{\left[\mathbf{P}_{h}\left(\mathbf{u}_{k}\right)\right]}\right)}{\prod\limits_{i=1}^{k}\left(\alpha_{i}-\alpha_{i+1}\right)!}\right]=\sum\limits_{\substack{ \alpha_{1}+\ldots+\alpha_{k}=k \\ \alpha_{1}\geq\ldots\geq\alpha_{k+1}=0 \\ \alpha_{i}\in \N\bigcup \left\{0\right\}}}\cfrac{1}{\prod\limits_{i=1}^{k}\left(\alpha_{i}-\alpha_{i+1}\right)!}\times
    \end{array}
\end{equation}
$$\times \left[\left\{\mathbf{N}^{\left(\alpha_{1}\right)}\left(\mathbf{u}_{0}\right)-\mathbf{N}^{\left(\alpha_{1}\right)}\left(\mathbf{P}_{h}\left(\mathbf{u}_{0}\right)\right)\right\}\left(\overset{\alpha_{1}-\alpha_{2}}{\left[\mathbf{u}_{1}\right]},\ldots, \overset{\alpha_{k}-\alpha_{k+1}}{\left[\mathbf{u}_{k}\right]}\right)+\right.$$
$$+\sum\limits_{i=1}^{\alpha_{1}-\alpha_{2}}\mathbf{N}^{\left(\alpha_{1}\right)}\left(\mathbf{P}_{h}\left(\mathbf{u}_{0}\right)\right)\left(\mathbf{u}_{1}-\mathbf{P}_{h}\left(\mathbf{u}_{1}\right),\overset{\alpha_{1}-\alpha_{2}-i}{\left[\mathbf{u}_{1}\right]}, \overset{i-1}{\left[\mathbf{P}_{h}\left(\mathbf{u}_{1}\right)\right]}, \overset{\alpha_{2}-\alpha_{3}}{\left[\mathbf{u}_{2}\right]},\ldots, \overset{\alpha_{k}-\alpha_{k+1}}{\left[\mathbf{u}_{k}\right]}\right)+$$
$$+\sum\limits_{i=1}^{\alpha_{2}-\alpha_{3}}\mathbf{N}^{\left(\alpha_{1}\right)}\left(\mathbf{P}_{h}\left(\mathbf{u}_{0}\right)\right)\left(\overset{\alpha_{1}-\alpha_{2}}{\left[\mathbf{P}_{h}\left(\mathbf{u}_{1}\right)\right]},\mathbf{u}_{2}-\mathbf{P}_{h}\left(\mathbf{u}_{2}\right),\overset{\alpha_{2}-\alpha_{3}-i}{\left[\mathbf{u}_{2}\right]}, \overset{i-1}{\left[\mathbf{P}_{h}\left(\mathbf{u}_{2}\right)\right]}, \ldots, \overset{\alpha_{k}-\alpha_{k+1}}{\left[\mathbf{u}_{k}\right]}\right)+$$
\begin{equation*}
\begin{split}
\ldots+\sum\limits_{i=1}^{\alpha_{k}}\mathbf{N}^{\left(\alpha_{1}\right)}\left(\mathbf{P}_{h}\left(\mathbf{u}_{0}\right)\right)\left(\overset{\alpha_{1}-\alpha_{2}}{\left[\mathbf{P}_{h}\left(\mathbf{u}_{1}\right)\right]},\ldots, \overset{\alpha_{k-1}-\alpha_{k}}{\left[\mathbf{P}_{h}\left(\mathbf{u}_{k-1}\right)\right]},\right.
\mathbf{u}_{k}-\mathbf{P}_{h}\left(\mathbf{u}_{k}\right), & \left.\left. \overset{\alpha_{k}-i}{\left[\mathbf{u}_{k}\right]}, \overset{i-1}{\left[\mathbf{P}_{h}\left(\mathbf{u}_{k}\right)\right]}\right)\right].
\end{split}
\end{equation*}
From condition (\ref{operat_lema_21}), the convexity of the set $\G_{n}$ and from properties (\ref{operat_lema_22}) -- (\ref{operat_lema_22''}) of the operator $\mathbf{P}_{h}\left(\mathbf{u}\right)$ we would obtain
\begin{equation}\label{operat_lema_28}
\begin{array}{c}
  \left\|\left\{\mathbf{N}^{\left(\alpha_{1}\right)}\left(\mathbf{u}_{0}\right)-\mathbf{N}^{\left(\alpha_{1}\right)}\left(\mathbf{P}_{h}\left(\mathbf{u}_{0}\right)\right)\right\}\left(\overset{\alpha_{1}-\alpha_{2}}{\left[\mathbf{u}_{1}\right]},\ldots, \overset{\alpha_{k}-\alpha_{k+1}}{\left[\mathbf{u}_{k}\right]}\right)\right\|^{(2)}= \\[1.2em]
  =\left\|\mathbf{N}^{\left(\alpha_{1}+1\right)}\left(\mathbf{P}_{h}\left(\mathbf{u}_{0}\right)-\theta\left(\mathbf{P}_{h}\left(\mathbf{u}_{0}\right)-\mathbf{u}_{0}\right)\right)\right\|\times\\[1.2em]
  \times\left\|\left(\mathbf{u}_{0}-\mathbf{P}_{h}\left(\mathbf{u}_{0}\right)\right)\right\|^{(1)}\prod\limits_{i=1}^{k}\left(\left\|\mathbf{u}_{1}\right\|^{(1)}\right)^{\alpha_{i}-\alpha_{i+1}}\leq \\[1.2em]
   \leq h\widetilde{N}^{(\alpha_{1}+1)}\left(\left\|\mathbf{P}_{h}\left(\mathbf{u}_{0}\right)-\theta\left(\mathbf{P}_{h}\left(\mathbf{u}_{0}\right)-\mathbf{u}_{0}\right)\right\|^{(1)}\right)\left\|\mathbf{u}_{0}\right\|^{(1)}_{1}\prod\limits_{i=1}^{k}\left(\left\|\mathbf{u}_{1}\right\|^{(1)}\right)^{\alpha_{i}-\alpha_{i+1}}\leq\\[1.2em]
   \leq h\widetilde{N}^{(\alpha_{1}+1)}\left(\left\|\mathbf{u}_{0}\right\|^{(1)}\right)\left\|\mathbf{u}_{0}\right\|^{(1)}_{1}\prod\limits_{i=1}^{k}\left(\left\|\mathbf{u}_{1}\right\|^{(1)}\right)^{\alpha_{i}-\alpha_{i+1}}\leq \\[1.2em]
   \leq h\widetilde{N}^{(\alpha_{1}+1)}\left(\ltri|\mathbf{u}_{0}\rtri|\right)\ltri|\mathbf{u}_{0}\rtri|\prod\limits_{i=1}^{k}\left(\ltri|\mathbf{u}_{1}\rtri|\right)^{\alpha_{i}-\alpha_{i+1}},\quad \theta\in\left(0, 1\right).
\end{array}
\end{equation}
Taking into account (\ref{operat_lema_27}) and (\ref{operat_lema_28}), we would  get the inequality
\begin{equation}\label{operat_lema_29}
\begin{array}{c}
  \left\|A_{k}\left(\mathbf{N}\left(\cdot\right); \left[\mathbf{u}_{i}\right]_{i=0}^{k}\right)-A_{k}\left(\mathbf{N}\left(\cdot\right); \left[\mathbf{P}_{h}\left(\mathbf{u}_{i}\right)\right]_{i=0}^{k}\right)\right\|^{(1)}\leq \\[1.2em]
  \leq\sum\limits_{\substack{ \alpha_{1}+\ldots+\alpha_{k}=k \\ \alpha_{1}\geq\ldots\geq\alpha_{k+1}=0 \\ \alpha_{i}\in \N\bigcup \left\{0\right\}}}\cfrac{1}{\prod\limits_{i=1}^{k}\left(\alpha_{i}-\alpha_{i+1}\right)!}\times\left[h\widetilde{N}^{(\alpha_{1}+1)}\left(\ltri|\mathbf{u}_{0}\rtri|\right)\ltri|\mathbf{u}_{0}\rtri|\prod\limits_{i=1}^{k}\left(\ltri|\mathbf{u}_{i}\rtri|\right)^{\alpha_{i}-\alpha_{i+1}}+\right.\\[1.2em]
\end{array}
\end{equation}
$$\begin{array}{c}
    +\left\|\mathbf{N}^{\left(\alpha_{1}\right)}\left(\mathbf{P}_{h}\left(\mathbf{u}_{0}\right)\right)\right\|\left\|\mathbf{u}_{1}-\mathbf{P}_{h}\left(\mathbf{u}_{1}\right)\right\|^{(1)}\prod\limits_{i=2}^{k} \left(\left\|\mathbf{u}_{i}\right\|^{(1)}\right)^{\alpha_{i}-\alpha_{i+1}}\times \\[1.2em]
    \times\left\{\sum\limits_{i=1}^{\alpha_{1}-\alpha_{2}}\left(\left\|\mathbf{u}_{1}\right\|^{(1)}\right)^{\alpha_{1}-\alpha_{2}-i} \left(\left\|\mathbf{P}_{h}\left(\mathbf{u}_{1}\right)\right\|^{(1)}\right)^{i-1}\right\}+\ldots + \\[1.2em]
    +\left\|\mathbf{N}^{\left(\alpha_{1}\right)}\left(\mathbf{P}_{h}\left(\mathbf{u}_{0}\right)\right)\right\|\left\|\mathbf{u}_{k}-\mathbf{P}_{h}\left(\mathbf{u}_{k}\right)\right\|^{(1)}\prod\limits_{i=1}^{k-1} \left(\left\|\mathbf{P}_{h}\left(\mathbf{u}_{i}\right)\right\|^{(1)}\right)^{\alpha_{i}-\alpha_{i+1}}\times \\[1.2em]
    \left.\times\left\{\sum\limits_{i=1}^{\alpha_{k}}\left(\left\|\mathbf{u}_{k}\right\|^{(1)}\right)^{\alpha_{k}-i} \left(\left\|\mathbf{P}_{h}\left(\mathbf{u}_{k}\right)\right\|^{(1)}\right)^{i-1}\right\}\right]\leq
  \end{array}$$
  $$\begin{array}{c}
      \leq h\sum\limits_{\substack{ \alpha_{1}+\ldots+\alpha_{k}=k \\ \alpha_{1}\geq\ldots\geq\alpha_{k+1}=0 \\ \alpha_{i}\in \N\bigcup \left\{0\right\}}}\cfrac{1}{\prod\limits_{i=1}^{k}\left(\alpha_{i}-\alpha_{i+1}\right)!}\prod\limits_{i=1}^{k}\left(\ltri|\mathbf{u}_{i}\rtri|\right)^{\alpha_{i}-\alpha_{i+1}}\times \\[1.2em]
      \times\left[\widetilde{N}^{(\alpha_{1}+1)}\left(\ltri|\mathbf{u}_{0}\rtri|\right)\ltri|\mathbf{u}_{0}\rtri|+\alpha_{1}\widetilde{N}^{(\alpha_{1})}\left(\ltri|\mathbf{u}_{0}\rtri|\right)\right]=hA_{k}\left(\widetilde{N}^{(1)}\left(\cdot\right)\times\left(\cdot\right); \left[\ltri|\mathbf{u}_{i}\rtri|\right]_{i=0}^{k}\right),\\[1.2em]
      \forall \mathbf{u}_{0}\in \G_{n}\bigcup \Eg_{1},\quad \forall \mathbf{u}_{i}\in \Eg_{1},\quad i\in \overline{1,k},
    \end{array}
  $$
which we needed to prove. The theorem is proved. $\blacksquare$

 To prove the main result of the paper we need to use the partial case of lemma  \ref{L_pro_Adom_1_Operat} which we are going to formulate below. For this reason we consider the sets $V_{m}\left(\mathbb{Q}\left[t_{0}, +\infty\right)\right)$ and $M_{m}\left(\mathbb{Q}\left[t_{0}, +\infty\right)\right).$ The former set is the set of column vectors of the dimension $m$ and the last one is the set of square matrices of the order $m\in\N.$ The components of the vectors from $V_{m}\left(\mathbb{Q}\left[t_{0}, +\infty\right)\right)$ and matrices from $M_{m}\left(\mathbb{Q}\left[t_{0}, +\infty\right)\right)$ are real valued piecewise continuous functions defined on the interval  $\left[t_{0}, +\infty\right).$ We will also refer to them as to the linear spaces $V_{m}\left(\mathbb{Q}\left[t_{0}, +\infty\right)\right)$ and $M_{m}\left(\mathbb{Q}\left[t_{0}, +\infty\right)\right),$ keeping in mind the sets of vectors/matrices  mentioned above together with the addition operation and the operation of multiplication of the vector/matrix by a real number  (see, for example, \cite[p. 13--15]{Gantmaher}). Let  $\left<\overrightarrow{u}, \overrightarrow{v}\right>$ denotes the scalar product of the vectors $\overrightarrow{u}, \overrightarrow{v}$ from  $V_{m}\left(\R\right).$
 We also define
 \begin{equation}\label{operat_lema_30}
 \begin{array}{c}
   \ltri|\overrightarrow{u}\left(t\right)\rtri|_{0, \left[a, b\right)}=\sup\limits_{t\in\left[a, b\right)}\left\|\overrightarrow{u}\left(t\right)\right\|=\sup\limits_{t\in\left[a, b\right)}\sqrt{\left<\overrightarrow{u}\left(t\right), \overrightarrow{u}\left(t\right)\right>}, \\[1.2em]
   \forall \overrightarrow{u}\left(t\right)\in V_{m}\left(\mathbb{Q}\left[t_{0}, +\infty\right)\right),
 \end{array}
 \end{equation}
 \label{P_funk_vect_norm}
 \begin{equation}\label{operat_lema_31}
 \begin{array}{c}
   \ltri|A\left(t\right)\rtri|_{0, \left[a, b\right)}=\sup\limits_{t\in\left[a, b\right)}\Bigg(\sup\limits_{\overrightarrow{u}\in V_{m}\left(\R\right)}\cfrac{\left\|A\left(t\right)\overrightarrow{u}\right\|}{\left\|\overrightarrow{u}\right\|}\Bigg), \\[1.2em]
   \forall A\left(t\right)\in M_{m}\left(\mathbb{Q}\left[t_{0}, +\infty\right)\right),\quad \varnothing\neq \left[a, b\right)\subseteq \left[t_{0}, +\infty\right).
 \end{array}
 \end{equation}

 It is easy to verify that the function $\ltri|\cdot\rtri|_{0, \left[a, b\right)},$ defined by formulas (\ref{operat_lema_30}), (\ref{operat_lema_31}) is the norm of the linear spaces $V_{m}\left(\mathbb{Q}\left[a, b\right)\right)$ and $M_{m}\left(\mathbb{Q}\left[a, b\right)\right).$ For convenience we will use the abbreviation
 \begin{equation}\label{operat_lema_32}
    \ltri|\cdot \rtri|_{0}\overset{\textmd{def}}{=}\ltri|\cdot\rtri|_{0, \left[a, b\right)},\quad \mbox{when} \quad \left[a, b\right)=\left[t_{0}, +\infty\right).
 \end{equation}

Let us fix a grid \eqref{Chapt_3_1_eq_4} on the interval $\left[t_{0}, +\infty\right)$
 and consider a set (a linear space) $V_{m}\left(\mathbb{Q}_{\omega}^{1}\left[t_{0}, +\infty\right)\right),$ which consists of the column vectors of dimension $m.$ The components of the vectors from $V_{m}\left(\mathbb{Q}_{\omega}^{1}\left[t_{0}, +\infty\right)\right)$ are continuous real valued functions defined on the interval $\left[t_{0}, +\infty\right)$ and continuously differentiable on each subinterval $\left(t_{i-1}, t_{i}\right),$ $\forall i\in\N,$ (we admit that there exist finite left and right derivatives at the points $t_{i}, i\in\N$ which can be unequal).
Then we can define
\begin{equation}\label{operat_lema_37}
    \ltri|\overrightarrow{u}\left(t\right)\rtri|_{1}\equiv \ltri|\cfrac{d}{d t}\overrightarrow{u}\left(t\right)\rtri|_{0},\quad \forall \overrightarrow{u}\left(t\right)\in V_{m}\left(\mathbb{Q}_{\omega}^{1}\left[t_{0}, +\infty\right)\right),
\end{equation}
 \label{P_funk_vect_norm_1}
\begin{equation}\label{operat_lema_38}
    \ltri|\overrightarrow{u}\left(t\right)\rtri|\equiv \max\left\{\ltri|\overrightarrow{u}\left(t\right)\rtri|_{0}, \ltri|\overrightarrow{u}\left(t\right)\rtri|_{1}\right\}, \quad \forall \overrightarrow{u}\left(t\right)\in V_{m}\left(\mathbb{Q}_{\omega}^{1}\left[t_{0}, +\infty\right)\right).
\end{equation}
 It is easy to verify that the functional $\ltri|\cdot\rtri|_{1}$ (\ref{operat_lema_37}) satisfies the conditions of seminorm in the space $V_{m}\left(\mathbb{Q}_{\omega}^{1}\left[t_{0}, +\infty\right)\right),$ thus, the functional $\ltri|\cdot\rtri|$ \eqref{operat_lema_38} is the norm in $V_{m}\left(\mathbb{Q}_{\omega}^{1}\left[t_{0}, +\infty\right)\right).$

We define the map $\mathbf{P}_{h}\colon V_{m}\left(\mathbb{Q}_{\omega}^{1}\left[t_{0}, +\infty\right)\right) \rightarrow V_{m}\left(\mathbb{Q}\left[t_{0}, +\infty\right)\right)$ in the following way
\begin{equation}\label{operat_lema_36}
\begin{array}{c}
  \mathbf{P}_{h}\left(\overrightarrow{u}\left(t\right)\right)=\overrightarrow{u}\left(t_{i-1}\right),\; \forall t\in\left[t_{i-1}, t_{i}\right),\; \forall i\in\N,\;\forall \overrightarrow{u}\left(t\right)\in V_{m}\left(\mathbb{Q}_{\omega}^{1}\left[t_{0}, +\infty\right)\right).
\end{array}
\end{equation}

Now we need to prove some auxiliary statements concerned with the operator $\mathbf{P}_{h}.$
\begin{lemma}\label{A_Lemma}
Let for some $\theta_{1}>0,\; \mathbf{u}, \mathbf{v}\in \E_{1}$ we have
\begin{equation}\label{L2.1}
  \left<\mathbf{u}+\theta_1\mathbf{v},\mathbf{u}+\theta_1\mathbf{v}\right>\ge\left<\mathbf{u},\mathbf{u}\right>,
\end{equation}
then for all $ \theta^{\prime},\theta''$ such that $\theta'\geq\theta''\geq\frac{1}{2}\theta_{1}$ the inequality
\begin{equation}\label{L2.2}
    \left<\mathbf{u}+\theta'\mathbf{v},\mathbf{u}+\theta'\mathbf{v}\right>\geq\left<\mathbf{u}+\theta''\mathbf{v},\mathbf{u}+\theta''\mathbf{v}\right>.
\end{equation}
holds true.
\end{lemma}
\emph{Proof.}
From (\ref{L2.1}) we would get
$$
\left<\mathbf{u},\mathbf{u}\right>+2\theta_1\Re\left(\left<\mathbf{u},\mathbf{v}\right>\right)+\theta_{1}^2\left<\mathbf{v},\mathbf{v}\right>\geq\left<\mathbf{u},\mathbf{u}\right>.
$$
From the last inequality, taking into account that $\theta_1>0$, we can obtain
\begin{equation}\label{L2.3}
   \Re\left(\left<\mathbf{u},\mathbf{v}\right>\right)\geq-\frac{\theta_1}{2}\left<\mathbf{v},\mathbf{v}\right>.
\end{equation}
Let us consider the first derivative of the function
$$
F\left(\theta\right)=\left<\mathbf{u}+\theta\mathbf{v},\mathbf{u}+\theta\mathbf{v}\right>
$$
with respect to $\theta:$

\begin{equation}\label{L2.4}
    F'\left(\theta\right)=\left<\mathbf{u}+\theta\mathbf{v},\mathbf{u}+\theta\mathbf{v}\right>_{\theta}^{'}=2\Re\left(\left<\mathbf{v},\mathbf{u}+\theta\mathbf{v}\right>\right)=2\Re\left(\left<\mathbf{v},\mathbf{u}\right>\right)+2\theta\left<\mathbf{v},\mathbf{v}\right>
\end{equation}
We are going to prove that $\forall\; \theta\geq\frac{1}{2}\theta_{1}$ the expression \eqref{L2.4} is nonnegative.
Indeed, using the inequality (\ref{L2.3}) from (\ref{L2.4}) we would get
$$
F'\left(\theta\right)\geq-\theta_1\left<\mathbf{v},\mathbf{v}\right>+2\theta\left<\mathbf{v},\mathbf{v}\right>=\left(2\theta-\theta_{1}\right)\left<\mathbf{v},\mathbf{v}\right>\geq 0.
$$
The last inequality means that on the interval $\left[\frac{1}{2}\theta_1,+ \infty\right)$ the function $F\left(\theta\right)$ is nondecreasing and  $\forall \; \theta',\theta'' \colon \theta'\ge\theta''\geq\frac{1}{2}\theta_{1}$ the inequality $F\left(\theta'\right)\geq F\left(\theta''\right)$ holds true. The lemma is proved. $\blacksquare$

\begin{lemma}\label{B_Lemma}
 For all $ \overrightarrow{u}(t)\in V_{m}\left(Q^{1}_{\omega}\left[t_{0}, +\infty\right)\right),$ for all $ \theta_{1},\theta_{2}$ such that $\theta_{1}\geq\theta_{2}$ and for the operator $\mathbf{P}_{h}$ defined in \eqref{operat_lema_36} the following inequality

\begin{equation}\label{Tv.5_1}
\ltri|\mathbf{P}_{h}\left(\overrightarrow{u}\left(t\right)\right)-\theta_{1}\left(\mathbf{P}_{h}\left(\overrightarrow{u}\left(t\right)\right)-\overrightarrow{u}\left(t\right)\right)\rtri|_{0}\geq
\ltri|\mathbf{P}_{h}\left(\overrightarrow{u}\left(t\right)\right)-\theta_{2}\left(\mathbf{P}_{h}\left(\overrightarrow{u}\left(t\right)\right)-\overrightarrow{u}\left(t\right)\right)\rtri|_{0}
\end{equation}
holds true.
\end{lemma}

{\bf Proof.}
Let us fix an arbitrary $n\in \N.$
As a first step we are going to prove the validity of the following inequality
\begin{equation}\label{Tv.5}
\begin{array}{c}
  \ltri|\mathbf{P}_{h}\left(\overrightarrow{u}\left(t\right)\right)-\theta_{1}\left(\mathbf{P}_{h}\left(\overrightarrow{u}\left(t\right)\right)-\overrightarrow{u}\left(t\right)\right)\rtri|_{0,[t_0,t_n)}\geq \\ [1.2em]
 \geq \ltri|\mathbf{P}_{h}\left(\overrightarrow{u}\left(t\right)\right)-\theta_{2}\left(\mathbf{P}_{h}\left(\overrightarrow{u}\left(t\right)\right)-\overrightarrow{u}\left(t\right)\right)\rtri|_{0,[t_0,t_n)}.
 \end{array}
\end{equation}
Let $t^{\ast}\in[t_0,t_n]$ be a point such that
$$\begin{array}{c}
    A=\ltri|\mathbf{P}_{h}\left(\overrightarrow{u}\left(t\right)\right)-\theta_{2}\left(\mathbf{P}_{h}\left(\overrightarrow{u}\left(t\right)\right)-\overrightarrow{u}\left(t\right)\right)\rtri|_{0,[t_0,t_n)}= \\[1.2em]
    =\left\|\mathbf{P}_{h}\left(\overrightarrow{u}\left(t^{\ast}\right)\right)-\theta_{2}\left(\mathbf{P}_{h}\left(\overrightarrow{u}\left(t^{\ast}\right)\right)-\overrightarrow{u}\left(t^{\ast}\right)\right)\right\|. \end{array}
$$
If $t^{\ast}=t_{n}$ then
$$A=\left\|\overrightarrow{u}\left(t_{n}\right)\right\|\leq \ltri|\mathbf{P}_{h}\left(\overrightarrow{u}\left(t\right)\right)-\theta_{1}\left(\mathbf{P}_{h}\left(\overrightarrow{u}\left(t\right)\right)-\overrightarrow{u}\left(t\right)\right)\rtri|_{0,[t_0,t_n)}\quad \forall\; \theta_{1} \in \R$$ and inequality \eqref{Tv.5} is proved.

Now we consider the case when $t^{\ast}\neq t_{n}.$ Without loss of generality we can assume that
 $t^{\ast}\in[t_0,t_1).$

It is easy to see that  $A\geq\left\|\overrightarrow{u}\left(t_{0}\right)\right\|,$ thus
\begin{equation}\label{Tv.5_2}
    \begin{array}{c}
      A=\left\|\mathbf{P}_{h}\left(\overrightarrow{u}\left(t^{\ast}\right)\right)-\theta_{2}\left(\mathbf{P}_{h}\left(\overrightarrow{u}\left(t^{\ast}\right)\right)-\overrightarrow{u}\left(t^{\ast}\right)\right)\right\|= \\  [1.2em]
      =\left\|\overrightarrow{u}\left(t_{0}\right)-\theta_{2}\left(\overrightarrow{u}\left(t_{0}\right)-\overrightarrow{u}\left(t^{\ast}\right)\right)\right\|\geq \left\|\overrightarrow{u}\left(t_{0}\right)\right\|.
    \end{array}
\end{equation}
Inequality \eqref{Tv.5_2} means that the assumptions of lemma \ref{A_Lemma} would be satisfied if we put $\E_{1} = V_{m}\left(Q^{1}_{\omega}\left[t_{0}, +\infty\right)\right),$ $\mathbf{u}=\overrightarrow{u}\left(t_{0}\right)$ and $\mathbf{v}=\overrightarrow{u}\left(t_{0}\right)-\overrightarrow{u}\left(t^{\ast}\right).$ Thus lemma \ref{A_Lemma} yields
$$
\begin{array}{c}
  A=\left\|\overrightarrow{u}\left(t_{0}\right)-\theta_{2}\left(\overrightarrow{u}\left(t_{0}\right)-\overrightarrow{u}\left(t^{\ast}\right)\right)\right\|\leq \left\|\overrightarrow{u}\left(t_{0}\right)-\theta_{1}\left(\overrightarrow{u}\left(t_{0}\right)-\overrightarrow{u}\left(t^{\ast}\right)\right)\right\| \leq \\[1.2em]
   \leq \ltri|\mathbf{P}_{h}\left(\overrightarrow{u}\left(t\right)\right)-\theta_{1}\left(\mathbf{P}_{h}\left(\overrightarrow{u}\left(t\right)\right)-\overrightarrow{u}\left(t\right)\right)\rtri|_{0,[t_0,t_1)}\leq
   \\[1.2em]
   \leq\ltri|\mathbf{P}_{h}\left(\overrightarrow{u}\left(t\right)\right)-\theta_{1}\left(\mathbf{P}_{h}\left(\overrightarrow{u}\left(t\right)\right)-\overrightarrow{u}\left(t\right)\right)\rtri|_{0,[t_0,t_n)}.
\end{array}
$$
So, inequality \eqref{Tv.5} has been proved.

Due to the freedom of choice of the constant $n\in \N$ inequality \eqref{Tv.5} yields \eqref{Tv.5_1}. It concludes the proof. $\blacksquare$

\begin{lemma}\label{L_pro_Adom_1}
    Let $\mathbf{N}\left(t, \overrightarrow{u}\right)\in M_{m}\left(C^{0, n}_{t, \overrightarrow{u}}\left(\left[t_{0}, +\infty\right)\times \R^{m}\right)\right)$ and there exists a scalar function $\widetilde{N}\left(u\right)\in C^{n}\left(\R\right)$ such that  $\forall \left(t, \overrightarrow{u}\right)\in \left[t_{0}, +\infty\right)\times \R^{m}$ the following inequalities hold true
   \begin{equation}\label{Chapt_3_2_eq_37}
     \sum\limits_{\substack{k_{1}+\ldots+k_{m}=k\\ k_{i}\in\N\bigcup \left\{0\right\}}}\cfrac{k !}{k_{1}!\ldots k_{m}!}\left\|\cfrac{\partial^{k}\mathbf{N}\left(t, \overrightarrow{u}\right)}{\partial u_{1}^{k_{1}}\ldots\partial u_{m}^{k_{m}}}\right\|\leq \left.\cfrac{d^{k}}{d u^{k}}\widetilde{N}\left(u\right)\right|_{u=\left\|\overrightarrow{u}\right\|},\;\forall k\in\overline{0, n},
   \end{equation}
   then
   \begin{equation}\label{Chapt_3_2_eq_38'}
   \begin{array}{c}
     \left\|A_{k}\left(\mathbf{N}\left(t,\cdot\right); \left[\overrightarrow{u}_{i}\left(t\right)\right]_{i=0}^{k}\right)\right\|\leq A_{k}\left(\widetilde{N}\left(\cdot\right); \left[\ltri|\overrightarrow{u}_{i}\left(t\right)\rtri|_{0}\right]_{i=0}^{k}\right), \\[1.2em]
     \forall \overrightarrow{u}_{i}\left(t\right)\in V_{m}\left(\mathbb{Q}\left[t_{0}, +\infty\right)\right),\quad i\in \overline{0,k},\quad \forall k\in \overline{0, n-1},
   \end{array}
   \end{equation}
   \begin{equation}\label{Chapt_3_2_eq_38}
   \begin{array}{c}
     \left\|A_{k}\left(\mathbf{N}\left(t, \cdot\right); \left[\overrightarrow{u}_{i}\left(t\right)\right]_{i=0}^{k}\right)-A_{k}\left(\mathbf{N}\left(t, \cdot\right); \left[\mathbf{P}_{h}\left(\overrightarrow{u}_{i}\left(t\right)\right)\right]_{i=0}^{k}\right)\right\|\leq \\ [1.2em]
     \leq h A_{k}\left(\widetilde{N}^{(1)}\left(\cdot\right)\times\left(\cdot\right); \ltri|\overrightarrow{u}_{0}\left(t\right)\rtri|,\ldots, \ltri|\overrightarrow{u}_{k}\left(t\right)\rtri|\right),\\[1.2em]
     \forall \overrightarrow{u}_{i}\left(t\right)\in V_{m}\left(\mathbb{Q}^{1}_{\omega}\left[t_{0}, +\infty\right)\right),\quad i\in \overline{0,k},\quad \forall k\in \overline{0, n-2}.
   \end{array}
   \end{equation}
\end{lemma}
{\bf Proof.} Let us show that the assertion of lemma \ref{L_pro_Adom_1} follows from lemma \ref{L_pro_Adom_1_Operat}. To achieve that we need to put $\E_{1}=\left<V_{m}\left(\mathbb{Q}\left[t_{0}, +\infty\right)\right); \ltri|\cdot \rtri|_{0}\right>,$ $\E_{2}=\left<M_{m}\left(\mathbb{Q}\left[t_{0}, +\infty\right)\right); \ltri|\cdot \rtri|_{0}\right>,$ $\G_{n}=V_{m}\left(\mathbb{Q}\left[t_{0}, +\infty\right)\right),$ $\mathcal{E}_{1}=\left<V_{m}\left(\mathbb{Q}^{1}_{\omega}\left[t_{0}, +\infty\right)\right); \ltri|\cdot\rtri|\right>.$ It is easy to see that the operator $\mathbf{P}_{h}$ (\ref{operat_lema_36}) satisfies the following inequality
\begin{equation}\label{operat_lema_39}
    \ltri|\overrightarrow{u}\left(t\right)-\mathbf{P}_{h}\left(\overrightarrow{u}\left(t\right)\right)\rtri|_{0}\leq h\ltri|\overrightarrow{u}\left(t\right)\rtri|_{1}, \quad \forall \overrightarrow{u}\left(t\right)\in\Eg_{1},
\end{equation}
moreover the inequality
\begin{equation}\label{operat_lema_41}
\begin{array}{c}
  \ltri|\mathbf{P}_{h}\left(\overrightarrow{u}\left(t\right)\right)-\theta_{1}\left(\mathbf{P}_{h}\left(\overrightarrow{u}\left(t\right)\right)-\overrightarrow{u}\left(t\right)\right)\rtri|_{0}\geq\ltri|\mathbf{P}_{h}\left(\overrightarrow{u}\left(t\right)\right)-\theta_{2}\left(\mathbf{P}_{h}\left(\overrightarrow{u}\left(t\right)\right)-\overrightarrow{u}\left(t\right)\right)\rtri|_{0}, \\[1.2em]
  \quad \forall \theta_{1}, \theta_{2}\in \R \colon \theta_{1}\geq \theta_{2}\geq 0,\quad\forall \overrightarrow{u}\left(t\right)\in \Eg_{1}
\end{array}
\end{equation}
is justified by lemma \ref{B_Lemma}.

 In the considered partial case inequalities (\ref{operat_lema_39}), (\ref{operat_lema_41})  are equivalent to inequalities (\ref{operat_lema_22}) and (\ref{operat_lema_22''}) respectively, and  additionally the condition similar to condition (\ref{operat_lema_22'}) is trivially satisfied. Let us show that inequalities (\ref{Chapt_3_2_eq_37}) implies (\ref{operat_lema_21}). Indeed, for arbitrary $\overrightarrow{u}\left(t\right), \overrightarrow{u}_{i}\left(t\right)\in \E_{1},$ $i\in \overline{1,k},$ $\forall t \in \left[t_{0}, +\infty\right)$ we will have (see. \cite[p. 264, p. 279]{shvartz_analiz})
\begin{equation}\label{operat_lema_42}
\begin{array}{c}
  \left\|\mathbf{N}^{(k)}_{\overrightarrow{u}}\left(t,\overrightarrow{u}\left(t\right)\right)\left(\overrightarrow{u}_{1}\left(t\right), \overrightarrow{u}_{2}\left(t\right), \ldots, \overrightarrow{u}_{k}\left(t\right)\right)\right\|= \\[1.2em]
  =\left\|\sum\limits_{i_{1}=1}^{m}u_{1, i_{1}}\left(t\right)\sum\limits_{i_{2}=1}^{m}u_{2, i_{2}}\left(t\right)\ldots \sum\limits_{i_{k}=1}^{m} u_{k, i_{k}}\left(t\right)\left.\cfrac{\partial^{k} \mathbf{N}\left(t,\overrightarrow{u}\right)}{\partial u_{i_{1}}\partial u_{i_{2}}\ldots \partial u_{i_{k}}}\right|_{\overrightarrow{u}=\overrightarrow{u}\left(t\right)}\right\|\leq \\[1.2em]
  \leq \sum\limits_{i_{1}=1}^{m}\left|u_{1, i_{1}}\left(t\right)\right|\ldots \sum\limits_{i_{k}=1}^{m} \left|u_{k, i_{k}}\left(t\right)\right|\left\|\left.\cfrac{\partial^{k} \mathbf{N}\left(t,\overrightarrow{u}\right)}{\partial u_{i_{1}}\partial u_{i_{2}}\ldots \partial u_{i_{k}}}\right|_{\overrightarrow{u}=\overrightarrow{u}\left(t\right)}\right\|\leq\\ [1.2em]
  \leq \left\|\overrightarrow{u}_{1}\left(t\right)\right\|\ldots \left\|\overrightarrow{u}_{k}\left(t\right)\right\|\sum\limits_{i_{1}=1}^{m}\ldots \sum\limits_{i_{k}=1}^{m} \left\|\left.\cfrac{\partial^{k} \mathbf{N}\left(t,\overrightarrow{u}\right)}{\partial u_{i_{1}}\partial u_{i_{2}}\ldots \partial u_{i_{k}}}\right|_{\overrightarrow{u}=\overrightarrow{u}\left(t\right)}\right\|=\\[1.2em]
  =\left(\prod\limits_{i=1}^{k}\left\|\overrightarrow{u}_{i}\left(t\right)\right\|\right)\sum\limits_{\substack{k_{1}+\ldots+k_{m}=k\\ k_{i}\in\N\bigcup \left\{0\right\}}}\cfrac{k !}{k_{1}!\ldots k_{m}!}\left\|\left.\cfrac{\partial^{k}\mathbf{N}\left(t, \overrightarrow{u}\right)}{\partial u_{1}^{k_{1}}\ldots\partial u_{m}^{k_{m}}}\right|_{\overrightarrow{u}=\overrightarrow{u}\left(t\right)}\right\|\leq \\[1.2em]
\end{array}
\end{equation}

\begin{equation*}
\begin{array}{c}
  \leq \left(\prod\limits_{i=1}^{k}\left\|\overrightarrow{u}_{i}\left(t\right)\right\|\right)\left.\cfrac{d^{k}}{d u^{k}}\widetilde{N}\left(u\right)\right|_{u=\left\|\overrightarrow{u}\left(t\right)\right\|}\leq \\[1.2em]
  \leq \left(\prod\limits_{i=1}^{k}\ltri|\overrightarrow{u}_{i}\left(t\right)\rtri|_{0}\right)\left.\cfrac{d^{k}}{d u^{k}}\widetilde{N}\left(u\right)\right|_{u=\ltri|\overrightarrow{u}\left(t\right)\rtri|_{0}},\quad \forall k\in \overline{0, n},
\end{array}
\end{equation*}
 where $\overrightarrow{u}_{i}\left(t\right)=\left[u_{i,1}\left(t\right),\dots, u_{i,m}\left(t\right)\right]^{T},$ $\forall i\in \overline{1,k}.$ Using definition  \ref{O_normi_klin_operatora} of the norm of $k$-linear operator, from (\ref{operat_lema_42}) it is easy to obtain inequalities (\ref{operat_lema_21}).

The presented above means that each condition of lemma \ref{L_pro_Adom_1_Operat} is satisfied, consequently its assertion, taking into account (\ref{operat_lema_37}) and (\ref{operat_lema_38}), implies inequalities (\ref{Chapt_3_2_eq_38'}) and (\ref{Chapt_3_2_eq_38}). The proof is complete. $\blacksquare$

\begin{conseq}\label{N_z lemi pro Adom1}
Let $\mathbf{N}\left(t, \overrightarrow{u}\right)\in M_{m}\left(C^{0, \infty}_{t, \overrightarrow{u}}\left(\left[t_{0}, +\infty\right)\times \R^{m}\right)\right)$ and there exists a scalar function $\widetilde{N}\left(u\right)\in C^{\infty}\left(\R\right)$ such that  $\forall \left(t, \overrightarrow{u}\right)\in \left[t_{0}, +\infty\right)\times \R^{m}$ the following inequalities hold true
   \begin{equation}\label{Chapt_3_2_eq_371}
     \sum\limits_{\substack{k_{1}+\ldots+k_{m}=k\\ k_{i}\in\N\bigcup \left\{0\right\}}}\cfrac{k !}{k_{1}!\ldots k_{m}!}\left\|\cfrac{\partial^{k}\mathbf{N}\left(t, \overrightarrow{u}\right)}{\partial u_{1}^{k_{1}}\ldots\partial u_{m}^{k_{m}}}\right\|\leq \left.\cfrac{d^{k}}{d u^{k}}\widetilde{N}\left(u\right)\right|_{u=\left\|\overrightarrow{u}\right\|},\;\forall k\in\N\bigcup \left\{0\right\}.
   \end{equation}
Then
\begin{equation}\label{Chapt_3_2_eq_51}
\begin{array}{c}
  \ltri|A_{k}\left(\mathbf{N}\left(t, \cdot\right); \left[\overrightarrow{u}_{i}\left(t\right)\right]_{i=0}^{k}\right)\rtri|_{0}\leq A_{k}\left(\widetilde{N}\left(\cdot\right); \left[\ltri|\overrightarrow{u}_{i}\left(t\right)\rtri|\right]_{i=0}^{k}\right), \\[1.2em]
  \forall t\in \left[t_{0}, +\infty\right),\quad \forall \overrightarrow{u}_{i}\left(t\right) \in V_{m}\left(\mathbb{Q}\left[t_{0}, +\infty\right)\right),\quad i\in\overline{0,k},\; \forall k\in\N\bigcup \left\{0\right\},
\end{array}
\end{equation}
\begin{equation}\label{Chapt_3_2_eq_43}
   \begin{array}{c}
     \ltri|A_{k}\left(\mathbf{N}\left(t, \cdot\right); \left[\overrightarrow{u}_{i}\left(t_{j-1}\right)\right]_{i=0}^{k}\right)-A_{k}\left(\mathbf{N}\left(t, \cdot\right); \left[\overrightarrow{u}_{i}\left(t\right)\right]_{i=0}^{k}\right)\rtri|_{0}\leq \\ [1.2em]
     \leq h A_{n}\left(\widetilde{N}^{(1)}\left(\cdot\right)\times \left(\cdot\right); \left[\ltri|\overrightarrow{u}_{i}\rtri|\right]_{i=0}^{k}\right),\quad \forall \overrightarrow{u}_{i}\left(t\right)\in V_{m}\left(\mathbb{Q}_{\omega}^{1}\left[t_{0}, +\infty\right)\right),\\ [1.2em]
     \forall t\in \left[t_{j-1}, t_{j}\right),\; j\in \N,\;i\in\overline{0,k},\; \forall k\in\N\bigcup \left\{0\right\}.
   \end{array}
\end{equation}
\end{conseq}

 The following lemma clarifies and generalizes lemma 2.2 from \cite{Gav_Laz_Mak_Sytn_2008}.
\begin{lemma}\label{L_pro Adom2}
For an arbitrary scalar function $\widetilde{N}\left(u\right)\in C^{\infty}\left(\R\right),$ $\forall u_{i}\in\R,$ $i\in\overline{0,n-1}$ the following equalities hold true
\begin{equation}\label{Chapt_3_2_eq_52}
\begin{array}{c}
  A_{n}\left(\widetilde{N}\left(\cdot\right); u_{0}, u_{1},\ldots, u_{n-1}, 0\right)= \\[1.2em]
  =\left.\cfrac{1}{n!}\cfrac{d^{n}}{d \tau^{n}}\left\{\widetilde{N}\left(\sum\limits_{i=0}^{\infty}\tau^{i}u_{i}\right)-\left.\cfrac{d \widetilde{N}\left(u\right)}{d u}\right|_{u=u_{0}}\sum\limits_{i=0}^{\infty}\tau^{i}u_{i}\right\}\right|_{\tau=0}, \quad \forall n\in \N.
\end{array}
\end{equation}
\end{lemma}

{\bf The convergence result for the FD-method applied to the Cauchy problem on the infinite interval.}
  We will use the notation $\mathbf{J}\left(t, \overrightarrow{u}\right)$ to describe the Jacobian matrix of a vector-valued function $\mathbf{N}\left(t, \overrightarrow{u}\right)\overrightarrow{u},$ $\mathbf{N}\left(t, \overrightarrow{u}\right)\in M_{m}\left(C^{0,1}_{t,\overrightarrow{u}}\left(\left[t_{0}, +\infty\right)\times \R^{m}\right)\right),$ dependent upon $\overrightarrow{u}$ (here we consider $t$ as a parameter), i.e.
  \begin{equation}\label{Chapt_3_2_eq_14}
  \begin{array}{c}
    \mathbf{J}\left(t, \overrightarrow{u}\right)=\mathbf{N}\left(t, \overrightarrow{u}\right)+\Bigg[\bigg(\cfrac{\partial}{\partial u_{1}}\mathbf{N}\left(t, \overrightarrow{u}\right)\bigg)\overrightarrow{u},\ldots, \bigg(\cfrac{\partial}{\partial u_{m}}\mathbf{N}\left(t, \overrightarrow{u}\right)\bigg)\overrightarrow{u}\Bigg], \\[1.2em]
    \overrightarrow{u}=\left[u_{1}, u_{2}, \ldots, u_{m}\right]^{T}.
  \end{array}
  \end{equation}

The following theorem gives the sufficient conditions for the convergence of FD-method applied to the Cauchy problem \eqref{Chapt_3_1_eq_1} on the interval $\left[t_{0}, +\infty\right)$ in the sense of definition \ref{O_FD_zbizhn}. In particular, this means that if we are going to apply the FD-method to the Cauchy problem \eqref{Chapt_3_1_eq_1}, which satisfies the conditions of the theorem stated below, then for an arbitrary interval $\left[t_{0}, t_{0}+H\right),$ $H>0$ we can use the uniform grid  $\widehat{\omega}$ \eqref{Chapt_3_1_eq_4} with the step  $h$ independent on $H$. It is easy to see that in such case the CPU time will be dependent on $H$ only linearly  (we assume that the computational complexity of the expressions for $\mathbf{N}\left(t, \overrightarrow{u}\right)$ and $\overrightarrow{\phi}\left(t\right)$ has the same order on  $\left[t_{0}, +\infty\right)$). The last fact is especially important when we are looking for the solution of the Cauchy problem \eqref{Chapt_3_1_eq_1} on a rather large interval.
\begin{theorem}\label{T_first_FD}
   Let the Cauchy problem \eqref{Chapt_3_1_eq_1} satisfies the following conditions
   \begin{enumerate}
   \item\label{U_1_FD_teor}  $\mathbf{N}\left(t, \overrightarrow{u}\right)=\sum\limits_{p=0}^{\infty}\sum\limits_{i_{1}+\ldots+i_{m}=p}u_{1}^{i_{1}}\ldots u_{m}^{i_{m}}\mathbf{N}_{i_{1}\ldots i_{m}}\left(t\right),\; \forall \left(t, \overrightarrow{u}\right)\in \left[t_{0}, +\infty\right)\times \R^{m},$
         where $\mathbf{N}_{i_{1}\ldots i_{m}}\left(t\right)\in M_{m}\left(C\left[t_{0}, +\infty\right)\right),$ $i_{k}\in \N\bigcup \left\{0\right\},$ $\forall k\in \overline{1,m},$ and in addition to that there exists a sequence of nonnegative real numbers $\left\{B_{i}\right\}_{i=0}^{\infty}$ such that $\sum\limits_{i_{1}+\ldots+i_{m}=p}\ltri|\mathbf{N}_{i_{1}\ldots i_{m}}\left(t\right)\rtri|_{0}\leq B_{p},$ $\forall p\in \N\bigcup \left\{0\right\},$ and the series $\sum\limits_{p=0}^{\infty}u^{p}B_{p}$ converges for all $u\in \R;$
    \item\label{U_2_FD_teor} $\overrightarrow{\phi}\left(t\right)\in V_{m}\left(C\left[t_{0}, +\infty\right)\right),$ $\ltri|\overrightarrow{\phi}\left(t\right)\rtri|_{0}\leq \kappa<+\infty;$
    \item\label{U_3_FD_teor} there exists a constant $\alpha\colon 0<\alpha\in\R$ such that $\forall \overrightarrow{v}, \overrightarrow{u}\in V_{m}\left(\R\right),$ $\forall t\in \left[t_{0}, +\infty\right)$ the inequality $$\left<\overrightarrow{v}, \mathbf{J}\left(t, \overrightarrow{u}\right)\overrightarrow{v}\right>\leq -\alpha \left<\overrightarrow{v}, \overrightarrow{v}\right> $$ is valid.
   \end{enumerate}
   Then for an arbitrary initial condition $\overrightarrow{u}_{0}\in V_{m}\left(\R\right)$ the solution of the Cauchy problem \eqref{Chapt_3_1_eq_1} exists and is unique on $\left[t_{0}, +\infty\right).$ The  FD-method for the Cauchy problem \eqref{Chapt_3_1_eq_1} converges to the exact solution of the problem with the following error estimates
   \begin{equation}\label{Chapt_3_2_eq_18}
    \ltri|\overrightarrow{u}\left(t\right)-\overset{p}{\overrightarrow{u}}\left(t\right)\rtri|_{0}\leq \frac{C}{\left(p+1\right)^{1+\varepsilon}}\frac{\left(h/R\right)^{p+1}}{1-h/R},\quad h<R,
   \end{equation}
   \begin{equation}\label{Chapt_3_2_eq_19}
    \ltri|\overrightarrow{u}\left(t\right)-\overset{p}{\overrightarrow{u}}\left(t\right)\rtri|_{0}\leq C\sum\limits_{j=p+1}^{\infty}\frac{1}{\left(j+1\right)^{1+\varepsilon}},\quad h=R,
   \end{equation}
    where the positive real constants  $C, R, \varepsilon$ depend on the input data of the problem \eqref{Chapt_3_1_eq_1} only.
\end{theorem}
To prove theorem \ref{T_first_FD} we need several auxiliary statements, which will be stated below.

\begin{lemma} \label{L_first_FD}
   Let $\mathbf{N}\left(t, \overrightarrow{u}\right)\in M_{m}\left(C^{0,1}_{t,\overrightarrow{u}}\left(\left[t_{0}, +\infty\right)\times \R^{m}\right)\right)$ and conditions \ref{U_2_FD_teor} and \ref{U_3_FD_teor} of theorem {\normalfont \ref{T_first_FD}} hold true,
then for an arbitrary constant $\varepsilon_{1}>0$  there exists a positive constant $\overline{h}=\overline{h}\left(\varepsilon_{1}\right)\in\R$ such that for every grid $\omega$ \eqref{Chapt_3_1_eq_1} with $h\leq \overline{h}$ the solution $\overrightarrow{u}^{(0)}\left(t\right)$ of the base problem \eqref{Chapt_3_1_eq_5} satisfies the inequality
  \begin{equation}\label{Chapt_3_2_eq_20}
  \begin{array}{c}
    \ltri|\overrightarrow{u}^{(0)}\left(t\right)\rtri|_{0}< \mu, \\ [1.2em]
    \mu=\max\left\{\left\|\overrightarrow{u}_{0}\right\|, \frac{\kappa}{\alpha}\right\}+\varepsilon_{1}.
  \end{array}
  \end{equation}
\end{lemma}
 {\bf Proof.} Let the assumptions of the lemma are satisfied. Let us set $i=1,$ fix an arbitrary  $t\in\left[t_{0}, t_{1}\right]$ and multiply both sides of the system (\ref{Chapt_3_1_eq_5}) on $\overrightarrow{u}^{(0)}\left(t\right).$ This results in the following (see \cite[p. 284 -- 288]{Demidovich})
\begin{equation}\label{Chapt_3_2_eq_21}
\begin{array}{c}
  \cfrac{1}{2}\cfrac{d}{d t}\left\|\overrightarrow{u}^{(0)}\left(t\right)\right\|^{2}=\cfrac{1}{2}\cfrac{d}{d t}\left<\overrightarrow{u}^{(0)}\left(t\right), \overrightarrow{u}^{(0)}\left(t\right)\right>= \\[1.2em]
  =\left<\overrightarrow{u}^{(0)}\left(t\right), \left(\mathbf{N}\left(t, \overrightarrow{u}^{(0)}\left(t\right)\right)\overrightarrow{u}^{(0)}\left(t\right)-\mathbf{N}\left(t, \overrightarrow{0}\right)\overrightarrow{0}\right)\right>+ \\[1.2em]
   +\left<\overrightarrow{u}^{(0)}\left(t\right), \left(\mathbf{N}\left(t, \overrightarrow{u}^{(0)}\left(t_{0}\right)\right)-\mathbf{N}\left(t, \overrightarrow{u}^{(0)}\left(t\right)\right)\right)\overrightarrow{u}^{(0)}\left(t\right)\right>+ \\ [1.2em]
\end{array}
\end{equation}
\begin{equation*}
\begin{array}{c}
   +\left<\overrightarrow{u}^{(0)}\left(t\right),\overrightarrow{\phi}\left(t\right)\right> = \left<\overrightarrow{u}^{(0)}\left(t\right), \mathbf{J}\left(t, \overrightarrow{u}_{1}\right)\overrightarrow{u}^{(0)}\left(t\right)\right>+\\[1.2em]
   +\left(t-t_{0}\right)\bigg<\overrightarrow{u}^{(0)}\left(t\right), \Upsilon\left(t,\overrightarrow{u}_{2},\overrightarrow{u}^{(0)}\left(t\right)\right)\bigg(\bigg.\cfrac{d}{d t}\overrightarrow{u}^{(0)}\left(t\right)\bigg|_{t=\overline{t}}\bigg)\bigg>+\left<\overrightarrow{u}^{(0)}\left(t\right),\overrightarrow{\phi}\left(t\right)\right>,\\
\end{array}
\end{equation*}
where $\overrightarrow{u}_{1}=\theta_{1}\overrightarrow{u}^{(0)}\left(t\right), $$\overrightarrow{u}_{2}=\overrightarrow{u}^{(0)}\left(t_{0}\right)+\theta_{2}\left(\overrightarrow{u}^{(0)}\left(t\right)-\overrightarrow{u}^{(0)}\left(t_{0}\right)\right), $ $\theta_{k}\in \left[0, 1\right],\;$ \mbox{$k=1,2,$} $\;\overline{t}\in \left[t_{0}, t\right],$
 $$\Upsilon\left(t, \overrightarrow{u}, \overrightarrow{v}\right)=\left[\cfrac{\partial \mathbf{N}\left(t, \overrightarrow{u}\right)}{\partial u_{1}}\overrightarrow{v},\;\cfrac{\partial \mathbf{N}\left(t, \overrightarrow{u}\right)}{\partial u_{2}}\overrightarrow{v},\ldots, \cfrac{\partial \mathbf{N}\left(t, \overrightarrow{u}\right)}{\partial u_{m}}\overrightarrow{v}\right].$$
 Taking the norms in the last equality of (\ref{Chapt_3_2_eq_21}) we get the following
\begin{equation}\label{Chapt_3_2_eq_22}
\begin{array}{c}
  \cfrac{1}{2}\cfrac{d}{d t}\left\|\overrightarrow{u}^{(0)}\left(t\right)\right\|^{2}\leq -\alpha \left\|\overrightarrow{u}^{(0)}\left(t\right)\right\|^{2}+\kappa\left\|\overrightarrow{u}^{(0)}\left(t\right)\right\|+ \\[1.2em]
  +h_{1}\left\|\overrightarrow{u}^{(0)}\left(t\right)\right\|\left\|\Upsilon\left(t,\overrightarrow{u}_{2},\overrightarrow{u}^{(0)}\left(t\right)\right)\right\|\times \\ [1.2em]
  \times \left\{\left\|\overrightarrow{u}^{(0)}\left(\overline{t}\right)\right\|\left\|\mathbf{N}\left(\overline{t}, \overrightarrow{u}\left(t_{0}\right)\right)\right\|+\kappa\right\}, \; h_{1}=t_{1}-t_{0}. \\ [1.2em]
\end{array}
\end{equation}

Let us fix an arbitrary $\varepsilon_{1}>0$ and, taking into account (\ref{Chapt_3_2_eq_20}), put
\begin{equation}\label{Chapt_3_2_eq_23}
 \mu_{1}=\max\left\{\left\|\overrightarrow{u}_{0}\right\|,\; \frac{\kappa}{\alpha}\right\}+\cfrac{\varepsilon_{1}}{2}.
\end{equation}
  It is easy to verify that in this case  $\alpha\mu_{1} -\kappa>0.$
We are going to prove that if
\begin{equation}\label{Chapt_3_2_eq_25}
    0<h_{1}\leq h\leq\overline{h}=\frac{\alpha\mu_{1} -\kappa}{\gamma_{2}\left(\mu\right)\left(\gamma_{1}\left(\mu\right)\mu+\kappa\right)},
\end{equation}
where \begin{equation}\label{Chapt_3_2_eq_24}
    \begin{array}{l}
      \gamma_{1}\left(\mu\right)=\sup\limits_{\overrightarrow{u}\in B_{\mu},\;
                     t\in\left[t_{0}, +\infty\right)
                   }\left\|\mathbf{N}\left(t,\overrightarrow{u}\right)\right\|<+\infty, \\[1.2em]
      \gamma_{2}\left(\mu\right)=\sup\limits_{\overrightarrow{u},\overrightarrow{v}\in B_{\mu},\;
                     t\in\left[t_{0}, +\infty\right)
      }\left\|\Upsilon\left(t,\overrightarrow{u},\overrightarrow{v}\right)\right\|<+\infty, \\[1.2em]
      B_{\mu}=\left\{\overrightarrow{u}\in V_{m}\left(\R\right) \mid \left\|\overrightarrow{u}\right\|\leq \mu\right\},
    \end{array}
\end{equation}
then
\begin{equation}\label{Chapt_3_2_eq_26}
    \sup\limits_{t\in \left[t_{0}, t_{1}\right]}\left\|\overrightarrow{u}^{(0)}\left(t\right)\right\|\leq \mu_{1}< \mu.
\end{equation}
 Indeed, let us assume that under condition (\ref{Chapt_3_2_eq_25}) inequality (\ref{Chapt_3_2_eq_26}) is not valid. So, there exists the real numbers  $t_{\ast},\; t^{\ast }\colon t_{0}<t_{\ast}< t^{\ast}\leq t_{1}$ such that
\begin{equation}\label{Chapt_3_2_eq_27}
    \mu_{1}<\left\|\overrightarrow{u}^{(0)}\left(t\right)\right\|< \mu,\; \forall t\in \left(t_{\ast}, t^{\ast}\right],\; \left\|\overrightarrow{u}^{(0)}\left(t_{\ast}\right)\right\|=\mu_{1}.
\end{equation}
Then from (\ref{Chapt_3_2_eq_22}) and (\ref{Chapt_3_2_eq_25}) $\forall t\in \left[t_{\ast}, t^{\ast}\right]$ we obtain
\begin{equation}\label{Chapt_3_2_eq_28}
    \frac{1}{2}\frac{d}{dt}\left\|\overrightarrow{u}^{(0)}\left(t\right)\right\|^{2}\leq \left\|\overrightarrow{u}^{(0)}\left(t\right)\right\|\left(-\alpha\mu_{1}+\kappa +\left(\alpha\mu_{1}-\kappa\right)\right)= 0.
\end{equation}
Inequality (\ref{Chapt_3_2_eq_28}) shows that the function $\left\|\overrightarrow{u}^{(0)}\left(t\right)\right\|$ is nonincreasing on $\left[t_{\ast}, t^{\ast}\right],$ hence $\left\|\overrightarrow{u}^{(0)}\left(t\right)\right\|\leq \mu_{1},$ $\forall t \in \left[t_{\ast}, t^{\ast}\right].$ The last inequality contradicts with assumption (\ref{Chapt_3_2_eq_27}). Thereby, inequality (\ref{Chapt_3_2_eq_26}) is valid.

Similarly, if we consider the base problem (\ref{Chapt_3_1_eq_5}) with $i=2,$ it is easy to prove, using \eqref{Chapt_3_2_eq_26}, that under condition (\ref{Chapt_3_2_eq_25}) the inequality
$$\sup\limits_{t\in \left[t_{1}, t_{2}\right]}\left\|\overrightarrow{u}^{(0)}\left(t\right)\right\|\leq \mu_{1}< \mu$$ holds true. In general, we can prove that condition (\ref{Chapt_3_2_eq_25}) together with the assumption $\left\|\overrightarrow{u}^{(0)}\left(t_{i-1}\right)\right\|\leq \mu_{1}$ implies the inequality
$$    \sup\limits_{t\in \left[t_{i-1}, t_{i}\right]}\left\|\overrightarrow{u}^{(0)}\left(t\right)\right\|\leq \mu_{1}< \mu.$$
Thus, using the principle of mathematical induction, we obtain
\begin{equation}\label{Chapt_3_2_eq_29}
    \sup\limits_{t\in \left[t_{0}, +\infty\right]}\left\|\overrightarrow{u}^{(0)}\left(t\right)\right\|\leq \mu_{1}< \mu.
\end{equation}
  From (\ref{Chapt_3_2_eq_29}) we obtain inequality (\ref{Chapt_3_2_eq_20}) for the $\overline{h}$ defined in (\ref{Chapt_3_2_eq_25}). The lemma is proved. $\blacksquare$

\begin{lemma}
  Let $\mathbf{J}\left(t\right)$ be the matrix from $M_{m}\left(C\left[t_{1}, t_{2}\right]\right),$ $t_{1}<t_{2}$ and
  \begin{equation}\label{Chapt_3_2_eq_30}
    \left<\overrightarrow{u}, \mathbf{J}\left(t\right)\overrightarrow{u}\right>\leq -\alpha\left<\overrightarrow{u},\overrightarrow{u}\right>,\; 0<\alpha\in \R,\quad \forall \overrightarrow{u}\in V_{m}\left(\R\right),\forall t\in\left[t_{1}, t_{2}\right],
  \end{equation}
  then the matrix $\Omega\left(t\right)\in M_{m}\left(C^{1}\left[t_{1}, t_{2}\right]\right)$ which is the solution of the Cauchy problem
  \begin{equation}\label{Chapt_3_2_eq_31}
    \frac{d}{d t}\Omega\left(t\right)-\mathbf{J}\left(t\right)\Omega\left(t\right)=0,\quad \Omega\left(t_{1}\right)=E, \quad t\in\left[t_{1}, t_{2}\right],
  \end{equation}
  satisfies the following inequalities
  \begin{equation}\label{Chapt_3_2_eq_32}
    \left\|\Omega\left(t\right)\right\|\leq e^{-\alpha \left(t-t_{1}\right)},\quad t\in \left[t_{1}, t_{2}\right],
  \end{equation}
  \begin{equation}\label{Chapt_3_2_eq_32'}
    \left\|\Omega^{-1}\left(t\right)\right\|\leq \exp\left(\left(t-t_{1}\right)\ltri|J\left(t\right)\rtri|_{0, \left[t_{1}, t_{2}\right)}\right),\quad t\in\left[t_{1}, t_{2}\right].
  \end{equation}

\end{lemma}

{\bf Proof.} Let the assumptions of the lemma hold true and $\Omega\left(t\right)$ denotes the solution of problem (\ref{Chapt_3_2_eq_31}), then $\forall \overrightarrow{u}\in V_{m}\left(\R\right)$ we have
\begin{equation}\label{Chapt_3_2_eq_33}
\begin{array}{c}
  \cfrac{1}{2}\cfrac{d}{d t}\left\|\Omega\left(t\right)\overrightarrow{u}\right\|^{2}=\cfrac{1}{2}\cfrac{d}{d t}\left<\Omega\left(t\right)\overrightarrow{u}, \Omega\left(t\right)\overrightarrow{u}\right>= \\[1.2em]
  =\left<\mathbf{J}\left(t\right)\Omega\left(t\right)\overrightarrow{u}, \Omega\left(t\right)\overrightarrow{u}\right>\leq -\alpha \left<\Omega\left(t\right)\overrightarrow{u}, \Omega\left(t\right)\overrightarrow{u}\right>=-\alpha \left\|\Omega\left(t\right)\overrightarrow{u}\right\|^{2}.
\end{array}
\end{equation}
From (\ref{Chapt_3_2_eq_33}) we obtain
\begin{equation}\label{Chapt_3_2_eq_34}
    \frac{\cfrac{d}{d t}\left\|\Omega\left(t\right)\overrightarrow{u}\right\|^{2}}{\left\|\Omega\left(t\right)\overrightarrow{u}\right\|^{2}}\leq -2\alpha.
\end{equation}
Integrating the both sides of inequality (\ref{Chapt_3_2_eq_34}) from $t_{1}$ to an arbitrary fixed point $t\in\left(t_{1}, t_{2}\right],$ we get the inequality
\begin{equation}\label{Chapt_3_2_eq_35}
    \ln\left(\left\|\Omega\left(t\right)\overrightarrow{u}\right\|^{2}\right)-\ln\left(\left\|\overrightarrow{u}\right\|^{2}\right)\leq -2\alpha\left(t-t_{1}\right).
\end{equation}
So far as inequality (\ref{Chapt_3_2_eq_35}) is valid for any arbitrary vector $\overrightarrow{u},$ we have
\begin{equation}\label{Chapt_3_2_eq_36}
\begin{array}{c}
  \sup\limits_{\overrightarrow{u}\in V_{m}\left(\R\right)}\ln\left(\cfrac{\left\|\Omega\left(t\right)\overrightarrow{u}\right\|^{2}}{\left\|\overrightarrow{u}\right\|^{2}}\right)=\ln\left(\sup\limits_{\overrightarrow{u}\in V_{m}\left(\R\right)}\cfrac{\left\|\Omega\left(t\right)\overrightarrow{u}\right\|^{2}}{\left\|\overrightarrow{u}\right\|^{2}}\right)= \\[1.2em]
  =\ln\left(\left\|\Omega\left(t\right)\right\|^{2}\right)\leq -2\alpha\left(t-t_{1}\right).
\end{array}
\end{equation}
Inequality (\ref{Chapt_3_2_eq_36}) implies (\ref{Chapt_3_2_eq_32}).

Now we are going to prove inequality (\ref{Chapt_3_2_eq_32'}). To do that we need to differentiate both parts of the identity $\Omega\left(t\right)\Omega^{-1}\left(t\right)=E,$ $\forall t\in\left[t_{1}, t_{2}\right]$ with respect to  $t.$ Then we get
\begin{equation}\label{Chapt_3_2_eq_32'_dodatok_1}
    \cfrac{d}{d t}\Omega^{-1}\left(t\right)=-\Omega^{-1}\left(t\right)J\left(t\right), \quad \Omega^{-1}\left(t_{1}\right)=E.
\end{equation}
Integrating the differential equation \eqref{Chapt_3_2_eq_32'_dodatok_1} from $t_{1}$ to an arbitrary point $t\in\left(t_{1}, t_{2}\right]$ we obtain
\begin{equation}\label{Chapt_3_2_eq_32'_dodatok_2}
    \Omega^{-1}\left(t\right)=E-\int\limits_{t_{1}}^{t}\Omega^{-1}\left(\xi\right)J\left(\xi\right)d\xi.
\end{equation}
Taking the norms in the last equality we come to the following inequality
\begin{equation}\label{Chapt_3_2_eq_32'_dodatok_3}
    \left\|\Omega^{-1}\left(t\right)\right\|\leq \left\|E\right\|+\int\limits_{t_{1}}^{t}\left\|\Omega^{-1}\left(\xi\right)\right\|\left\|J\left(\xi\right)\right\|d\xi=1+\int\limits_{t_{1}}^{t}\left\|\Omega^{-1}\left(\xi\right)\right\|\left\|J\left(\xi\right)\right\|d\xi.
\end{equation}
Applying the Gronwall–Bellman inequality (see \cite[p. 108]{Demidovich}) to  \eqref{Chapt_3_2_eq_32'_dodatok_3}   we obtain \eqref{Chapt_3_2_eq_32'}. The lemma is proved. $\blacksquare$

{\bf Proof (of the theorem \ref{T_first_FD}).} Let the conditions of the theorem hold true. The existence and uniqueness of the solution of the Cauchy problem (\ref{Chapt_3_1_eq_1}) on $\left[t_{0}, +\infty\right)$ was proved in \cite[p. 286]{Demidovich}.

Let us fix an infinite grid $\widehat{\omega}$ (\ref{Chapt_3_1_eq_4}). Then using a scalar parameter  $\tau\in\left[0, 1\right]$ we consider the following generalization of problem (\ref{Chapt_3_1_eq_1})
\begin{equation}\label{Chapt_3_2_eq_53}
\begin{array}{c}
  \cfrac{d}{d t}\overrightarrow{u}\left(t, \tau\right)-\mathbf{N}\left(t, \overrightarrow{u}\left( t_{i-1}, \tau\right)\right)\overrightarrow{u}\left(t, \tau\right)- \\[1.2em]
  -\tau\left\{\mathbf{N}\left(t, \overrightarrow{u}\left(t, \tau\right)\right)-\mathbf{N}\left(t, \overrightarrow{u}\left(t_{i-1},\tau\right)\right)\right\}\overrightarrow{u}\left(t, \tau\right)=\overrightarrow{\phi}\left(t\right),  \\[1.2em]
  \tau\in\left[0, 1\right],\; t\in \left[t_{i-1}, t_{i}\right),\; \forall i\in \N,
\end{array}
\end{equation}
\begin{equation}\label{Chapt_3_2_eq_54}
    \overrightarrow{u}\left(t_{0}, \tau\right)=\overrightarrow{u}_{0},\; \left[\overrightarrow{u}\left(t, \tau\right)\right]_{t=t_{i}}=0,\quad i\in \N, \quad \forall \tau\in \left[0, 1\right].
\end{equation}
It is easy to see that if $\tau=1$ then problem \eqref{Chapt_3_2_eq_53}, \eqref{Chapt_3_2_eq_54} transforms to \eqref{Chapt_3_1_eq_1}.

We assume that problem (\ref{Chapt_3_2_eq_53}), (\ref{Chapt_3_2_eq_54}) has a unique solution $\overrightarrow{u}\left(\tau, t\right),$ which can be expressed in the series form
\begin{equation}\label{Chapt_3_2_eq_55}
    \overrightarrow{u}\left(\tau, t\right)=\sum\limits_{i=0}^{\infty}\tau^{i} \overrightarrow{u}^{(i)}\left(t\right),\quad \forall \tau\in \left[0, 1\right], \quad \forall t\in\left[t_{0}, +\infty\right),
\end{equation}
where $\overrightarrow{u}^{(i)}\left(t\right)\in V_{m}\left(\mathbb{Q}_{\omega}^{1}\left[t_{0}, +\infty\right)\right), \; i\in\N\bigcup \left\{0\right\}$ (see notations on p. \pageref{P_funk_vect_norm}), furthermore
\begin{equation}\label{Chapt_3_2_eq_56}
\begin{split}
  \cfrac{\partial}{\partial t}\overrightarrow{u}\left(\tau, t\right)=\sum\limits_{i=0}^{\infty}\tau^{i} \cfrac{d}{d t}\overrightarrow{u}^{(i)}\left(t\right),\quad \forall \tau\in \left[0, 1\right], \quad \forall t\in \bigcup\limits_{i=1}^{\infty}\left(t_{i-1},\; t_{i}\right), \\
  \left.\cfrac{\partial}{\partial t}\overrightarrow{u}\left(\tau, t\right)\right|_{t=t_{j}}=\sum\limits_{i=0}^{\infty}\tau^{i} \lim\limits_{t\rightarrow t_{j}+0}\cfrac{d}{d t}\overrightarrow{u}^{(i)}\left(t\right),\quad \forall \tau\in\left[0, 1\right] \quad \forall j\in \N\bigcup \left\{0\right\}.
\end{split}
\end{equation}
Taking into account assumptions \eqref{Chapt_3_2_eq_55}, \eqref{Chapt_3_2_eq_56} and putting $\tau=0$ in \eqref{Chapt_3_2_eq_53}, \eqref{Chapt_3_2_eq_54}, we obtain the base problem \eqref{Chapt_3_1_eq_5} to define unknown term $\overrightarrow{u}^{(0)}\left(t\right).$ Similarly, if we substitute representation \eqref{Chapt_3_2_eq_55} into \eqref{Chapt_3_2_eq_53}, differentiate the obtained equality with respect to $\tau$ $j$ times ($j\in \N$) then divide both sides of the equality by $j!$ and finally put $\tau=0,$ we obtain the recursive system of Cauchy problems \eqref{Chapt_3_1_eq_6} to define unknown terms $\overrightarrow{u}^{(j)}\left(t\right),$ $\forall j\in\N.$

Let us prove that the FD-method for the Cauchy problem \eqref{Chapt_3_1_eq_1} converges in the sense of definition \ref{O_FD_zbizhn}. To do this we have to rewrite  \eqref{Chapt_3_1_eq_6} in the equivalent form, namely:
\begin{equation}\label{Chapt_3_2_eq_57}
\begin{array}{c}
  \cfrac{d}{d t}\overrightarrow{u}^{(j)}\left(t\right)-\left(\mathbf{N}\left(t, \overrightarrow{u}^{(0)}\left(t_{i-1}\right)\right)+\Upsilon_{i}\left(t, \overrightarrow{u}^{(0)}\left(t_{i-1}\right)\right)\right)\overrightarrow{u}^{(j)}\left(t\right)= \\[1.2em]
  =\Upsilon_{i}\left(t, \overrightarrow{u}^{(0)}\left(t\right)\right)\overrightarrow{u}^{(j)}\left(t_{i-1}\right)-\Upsilon_{i}\left(t, \overrightarrow{u}^{(0)}\left(t_{i-1}\right)\right)\overrightarrow{u}^{(j)}\left(t\right)+F^{(j)}\left(t\right),
\end{array}
\end{equation}
$ t\in\left[t_{i-1}, t_{i}\right),\; i,j\in\N,$ where, $\forall \overrightarrow{v}\in V_{m}\left(\R\right)$
\begin{equation}\label{Chapt_3_2_eq_14'}
    \Upsilon_{i}\left(t,\overrightarrow{v}\right)= \left.\Bigg[\bigg(\cfrac{\partial}{\partial u_{1}}\mathbf{N}\left(t, \overrightarrow{u}\right)\bigg)\overrightarrow{v},\ldots, \bigg(\cfrac{\partial}{\partial u_{m}}\mathbf{N}\left(t, \overrightarrow{u}\right)\bigg)\overrightarrow{v}\Bigg]\right|_{\overrightarrow{u}=\overrightarrow{u}^{(0)}\left(t_{i-1}\right)}.
\end{equation}
Taking into account notations \eqref{Chapt_3_2_eq_14}, we can rewrite equation \eqref{Chapt_3_2_eq_57} in the following form
\begin{equation}\label{Chapt_3_2_eq_58}
\begin{array}{c}
  \cfrac{d}{d t}\overrightarrow{u}^{(j)}\left(t\right)-\mathbf{J}_{i}\left(t\right)\overrightarrow{u}^{(j)}\left(t\right)=\Upsilon_{i}\bigg(t, \int\limits_{t_{i-1}}^{t}\cfrac{d}{d \xi}\overrightarrow{u}^{(0)}\left(\xi\right)d\xi\bigg)\overrightarrow{u}^{(j)}\left(t_{i-1}\right)- \\[1.2em]
  -\Upsilon_{i}\left(t, \overrightarrow{u}^{(0)}\left(t_{i-1}\right)\right)\int\limits_{t_{i-1}}^{t}\cfrac{d}{d \xi}\overrightarrow{u}^{(j)}\left(\xi\right)d \xi+F^{(j)}\left(t\right),\\[1.2em]
  \mathbf{J}_{i}\left(t\right)=\mathbf{J}\left(t, \overrightarrow{u}^{(0)}\left(t_{i-1}\right)\right),\quad t\in\left[t_{i-1}, t_{i}\right),\; i, j\in \N.
\end{array}
\end{equation}
 Then if we add to equation \eqref{Chapt_3_2_eq_58} the initial and matching conditions
\begin{equation}\label{Chapt_3_2_eq_59}
    \overrightarrow{u}^{(j)}\left(t_{0}\right)=\overrightarrow{u}_{0},\quad \left[\overrightarrow{u}^{(j)}\left(t\right)\right]_{t=t_{i}}=0,\; i,j\in\N,
\end{equation}
we would get the recursive system of the Cauchy problems with respect to the unknown vector-functions $\overrightarrow{u}^{(j)}\left(t\right),$ which is equivalent to system \eqref{Chapt_3_1_eq_6}.

Let us fix an arbitrary $j\in \N.$ Then the solution $\overrightarrow{u}^{(j)}\left(t\right)$ of the Cauchy problem \eqref{Chapt_3_2_eq_58}, \eqref{Chapt_3_2_eq_59} can be expressed in the form
\begin{equation}\label{Chapt_3_2_eq_60}
\begin{array}{c}
  \overrightarrow{u}^{(j)}\left(t\right)=\Bigg[\Omega_{i}\left(t\right)+\int\limits_{t_{i-1}}^{t}K_{i}\left(t,\xi\right)\Upsilon_{i}\bigg(\xi, \int\limits_{t_{i-1}}^{\xi}\cfrac{d}{d \eta}\overrightarrow{u}^{(0)}\left(\eta\right) d\eta\bigg)d\xi\Bigg]\overrightarrow{u}^{(j)}\left(t_{i-1}\right)- \\[1.2em]
  -\int\limits_{t_{i-1}}^{t}K_{i}\left(t, \xi\right)\Upsilon_{i}\left(\xi, \overrightarrow{u}^{(0)}\left(t_{i-1}\right)\right)\int\limits_{t_{i-1}}^{\xi}\cfrac{d}{d \eta}\overrightarrow{u}^{(j)}\left(\eta\right)d\eta d\xi+\\[1.2em]
  +\int\limits_{t_{i-1}}^{t}K_{i}\left(t, \xi\right)F^{(j)}\left(\xi\right)d\xi,\quad  t\in\left[t_{i-1}, t_{i}\right], \;i\in\N,
\end{array}
\end{equation}
where the matrices  $\Omega_{i}\left(t\right)\in M_{m}\left(C^{1}\left[t_{i-1}, t_{i}\right]\right)$ are the solutions of the Cauchy problems
$$\cfrac{d}{d t}\Omega_{i}\left(t\right)-\mathbf{J}_{i}\left(t\right)\Omega_{i}\left(t\right)=0,\quad \Omega_{i}\left(t_{i-1}\right)=E,\quad i\in \N,$$
$K_{i}\left(t, \xi\right)=\Omega_{i}\left(t\right) \Omega_{i}^{-1}\left(\xi\right).$

Lets use lemma \ref{L_first_FD} and assume that for some fixed constant $\varepsilon_{1}>0,$ a maximum step $h$ of the grid $\widehat{\omega}$ (\ref{Chapt_3_1_eq_4}) satisfies the inequality $h\leq\overline{h},$ where $\overline{h}=\overline{h}\left(\varepsilon_{1}\right)$ is a constant mentioned in lemma \ref{L_first_FD}. Since that, accordingly to the assertion of lemma \ref{L_first_FD}, we can assume that for the vector-function $\overrightarrow{u}^{(0)}\left(t\right)$ that is the solution of the base problem \eqref{Chapt_3_1_eq_5} inequality \eqref{Chapt_3_2_eq_20} is valid. Whereas the vector function $\overrightarrow{u}^{(0)}\left(t\right)$ depends on the grid $\widehat{\omega}$ (\ref{Chapt_3_1_eq_4}) it is worth to emphasize that the right side of estimate \eqref{Chapt_3_2_eq_20} does not depend on $\widehat{\omega}.$ After that using the following notations
\begin{equation}\label{Chapt_3_2_eq_60'}
    \begin{array}{c}
  B=\max\limits_{\substack{\left\|\overrightarrow{u}\right\|\leq \mu \\ \left\|\overrightarrow{v}\right\|\leq \mu}}\left\{\sum\limits_{i=1}^{m}\ltri|\cfrac{\partial \mathbf{N}\left(t, \overrightarrow{u}\right)}{\partial u_{i}}\rtri|_{0}\ltri|\mathbf{N}\left(t, \overrightarrow{u}\right)\overrightarrow{v}+\overrightarrow{\phi}\left(t\right)\rtri|_{0}\right\}<+\infty, \\
  C=\mu \max\limits_{\left\|\overrightarrow{u}\right\|\leq \mu}\left\{\sum\limits_{i=1}^{m}\ltri|\cfrac{\partial \mathbf{N}\left(t, \overrightarrow{u}\right)}{\partial u_{i}}\rtri|_{0}\right\}<+\infty
\end{array}
\end{equation}
and taking into account the evident inequality
$$\left\|\Upsilon_{i}\left(t, \overrightarrow{v}\left(t\right)\right)\overrightarrow{w}\left(t\right)\right\|=\left\|\bigg(\sum\limits_{i=1}^{m} w_{i}\left(t\right)\bigg.\cfrac{\partial \mathbf{N}\left(t, \overrightarrow{u}\right)}{\partial u_{i}}\bigg|_{\overrightarrow{u}=\overrightarrow{u}^{(0)}\left(t_{i-1}\right)}\bigg)\overrightarrow{v}\left(t\right)\right\|\leq$$
$$\leq \sum\limits_{i=1}^{m}\ltri|\bigg.\cfrac{\partial \mathbf{N}\left(t, \overrightarrow{u}\right)}{\partial u_{i}}\bigg|_{\overrightarrow{u}=\overrightarrow{u}^{(0)}\left(t_{i-1}\right)}\rtri|_{0}\ltri|\overrightarrow{v}\left(t\right)\rtri|_{0}\ltri|\overrightarrow{w}\left(t\right)\rtri|_{0},\quad \forall i\in \N,$$
$\forall \overrightarrow{v}\left(t\right), \overrightarrow{w}\left(t\right)\in V_{m}\left(C\left[t_{0}, +\infty\right)\right)$ from \eqref{Chapt_3_2_eq_60} we would obtain
\begin{equation}\label{Chapt_3_2_eq_61}
\begin{array}{c}
  \ltri|\overrightarrow{u}^{(j)}\left(t\right)\rtri|_{0, \left[t_{i-1}, t_{i}\right)}\leq\\[1.2em]
   \leq\left[1+\cfrac{h_{i}^{2}}{2}\sum\limits_{i=1}^{m}\ltri|\left.\cfrac{\partial \mathbf{N}\left(t, \overrightarrow{u}\right)}{\partial u_{i}}\right|_{\overrightarrow{u}=\overrightarrow{u}^{(0)}\left(t_{i-1}\right)}\rtri|_{0}\ltri|\cfrac{d}{d t}\overrightarrow{u}^{(0)}\left(t\right)\rtri|_{0} \right]\left\|\overrightarrow{u}^{(j)}\left(t_{i-1}\right)\right\|+\\[1.2em]
\end{array}
\end{equation}
\begin{equation*}
    \begin{array}{c}
  +\cfrac{h_{i}^{2}}{2}\sum\limits_{i=1}^{m}\ltri|\left.\cfrac{\partial \mathbf{N}\left(t, \overrightarrow{u}\right)}{\partial u_{i}}\right|_{\overrightarrow{u}=\overrightarrow{u}^{(0)}\left(t_{i-1}\right)}\rtri|_{0}\ltri|\overrightarrow{u}^{(0)}\left(t\right)\rtri|_{0}\ltri|\cfrac{d}{d t}\overrightarrow{u}^{(j)}\left(t\right)\rtri|_{0, \left[t_{i-1}, t_{i}\right)}+\\[1.2em]
  +h_{i}\ltri|F^{(j)}\left(t\right)\rtri|_{0}\leq
    \end{array}
\end{equation*}

$$\begin{array}{c}
    \leq \left[1+\cfrac{h_{i}^{2}B}{2}\right]\left\|\overrightarrow{u}^{(j)}\left(t_{i-1}\right)\right\|+\cfrac{h_{i}^{2}C}{2}\ltri|\cfrac{d}{d t}\overrightarrow{u}^{(j)}\left(t\right)\rtri|_{0, \left[t_{i-1}, t_{i}\right)}+h_{i}\ltri|F^{(j)}\left(t\right)\rtri|_{0}.
  \end{array}
$$
Similarly, from \eqref{Chapt_3_2_eq_60} we obtain
\begin{equation}\label{Chapt_3_2_eq_61'}
\begin{array}{c}
  \left\|\overrightarrow{u}^{(j)}\left(t_{i}\right)\right\|\leq \left[e^{-h_{i}\alpha}+\cfrac{h_{i}^{2}B}{2}\right]\left\|\overrightarrow{u}^{(j)}\left(t_{i-1}\right)\right\|+h_{i}\ltri|F^{(j)}\left(t\right)\rtri|_{0}+\\ [1.2em]
  +\cfrac{h_{i}^{2}C}{2}\ltri|\cfrac{d}{d t}\overrightarrow{u}^{(j)}\left(t\right)\rtri|_{0, \left[t_{i-1}, t_{i}\right)}, \quad i\in \N.
\end{array}
\end{equation}

On the other hand, from \eqref{Chapt_3_1_eq_6} we find
\begin{equation}\label{Chapt_3_2_eq_62}
\begin{array}{c}
  \overrightarrow{u}^{(j)}\left(t\right)=U_{i}\left(t\right)\overrightarrow{u}^{(j)}\left(t_{i-1}\right)+U_{i}\left(t\right)\int\limits_{t_{i-1}}^{t}U^{-1}_{i}\left(\xi\right)F^{(j)}\left(\xi\right)+ \\[1.2em]
  +U_{i}\left(t\right)\int\limits_{t_{i-1}}^{t}U^{-1}_{i}\left(\xi\right)\Upsilon_{i}\left(\xi, \overrightarrow{u}^{(0)}\left(\xi\right)\right)d\xi \overrightarrow{u}^{(j)}\left(t_{i-1}\right),
\end{array}
\end{equation}
where the matrix
$U_{i}\left(t\right)\in M_{m}\left(C^{1}\left[t_{i-1}, t_{i}\right]\right)$ is the solution of the Cauchy problem
$$\cfrac{d}{d t}U_{i}\left(t\right)-\mathbf{N}\left(t, \overrightarrow{u}^{(0)}\left(t_{i-1}\right)\right)\left(t\right)U_{i}\left(t\right)=0,\quad U_{i}\left(t_{i-1}\right)=E.$$

If we differentiate with respect to $t$ both sides of equality \eqref{Chapt_3_2_eq_62} we would obtain
\begin{equation}\label{Chapt_3_2_eq_63}
\begin{array}{c}
  \cfrac{d}{d t}\overrightarrow{u}^{(j)}\left(t\right)=\bigg[N_{i}\left(t\right)U_{i}\left(t\right)+\Upsilon_{i}\left(t, \overrightarrow{u}^{(0)}\left(t\right)\right)+\bigg. \\[1.2em]
  \bigg.+N_{i}\left(t\right)U_{i}\left(t\right)\int\limits_{t_{i-1}}^{t}U^{-1}_{i}\left(\xi\right)\Upsilon_{i}\left(\xi, \overrightarrow{u}^{(0)}\left(\xi\right)\right)d\xi\bigg]\overrightarrow{u}^{(j)}\left(t_{i-1}\right)+ \\[1.2em]
  +N_{i}\left(t\right)U_{i}\left(t\right)\int\limits_{t_{i-1}}^{t}U^{-1}_{i}\left(\xi\right)F^{(j)}\left(\xi\right)d\xi+F^{(j)}\left(t\right),\; t\in\left[t_{i-1}, t_{i}\right],\; i\in \N,
\end{array}
\end{equation}
 where $N_{i}\left(t\right)=\mathbf{N}\left(t, \overrightarrow{u}^{(0)}\left(t_{i-1}\right)\right).$ From \eqref{Chapt_3_2_eq_63} we get
\begin{equation}\label{Chapt_3_2_eq_64}
    \begin{array}{c}
      \ltri|\cfrac{d}{d t}\overrightarrow{u}^{(j)}\left(t\right)\rtri|_{0, \left[t_{i-1}, t_{i}\right)}\leq \Bigg[\ltri|N_{i}\left(t\right)\rtri|_{0}\exp\left\{h\ltri|N_{i}\left(t\right)\rtri|_{0}\right\}+\bigg.\\[1.2em]
      +\sum\limits_{i=1}^{m}\ltri|\Bigg.\cfrac{\partial \mathbf{N}\left(t, \overrightarrow{u}\right)}{\partial u_{i}}\bigg|_{\overrightarrow{u}=\overrightarrow{u}^{(0)}\left(t_{i-1}\right)}\rtri|_{0}\ltri|\overrightarrow{u}\left(t\right)\rtri|_{0}+h\ltri|N_{i}\left(t\right)\rtri|_{0}\exp\left\{2h\ltri|N_{i}\left(t\right)\rtri|_{0}\right\}\times \\[1.2em]
      \Bigg.\times\sum\limits_{i=1}^{m}\ltri|\bigg.\cfrac{\partial \mathbf{N}\left(t, \overrightarrow{u}\right)}{\partial u_{i}}\bigg|_{\overrightarrow{u}=\overrightarrow{u}^{(0)}\left(t_{i-1}\right)}\rtri|_{0}\ltri|\overrightarrow{u}\left(t\right)\rtri|_{0}\Bigg]\left\|\overrightarrow{u}^{(j)}\left(t_{i-1}\right)\right\|+\\[1.2em]
      +\left(h\ltri|N_{i}\left(t\right)\rtri|_{0}\exp\left\{2h\ltri|N_{i}\left(t\right)\rtri|_{0}\right\}+1\right)\ltri|F^{(j)}\left(t\right)\rtri|_{0}\leq\\[1.2em] \leq P\left\|\overrightarrow{u}^{(j)}\left(t_{i-1}\right)\right\|+Q\ltri|F^{(j)}\left(t\right)\rtri|_{0},
    \end{array}
\end{equation}
where
\begin{equation}\label{Chapt_3_2_eq_64'}
    P=Ne^{\overline{h}N}+CQ,\;Q=\overline{h}Ne^{2\overline{h}N}+1,\quad N=\max\limits_{\left\|\overrightarrow{u}\right\|\leq \mu}\ltri|\mathbf{N}\left(t, \overrightarrow{u}\right)\rtri|_{0}.
\end{equation}

From inequalities \eqref{Chapt_3_2_eq_61}, \eqref{Chapt_3_2_eq_61'} by virtue of \eqref{Chapt_3_2_eq_64} we obtain
\begin{equation}\label{Chapt_3_2_eq_65}
    \ltri|\overrightarrow{u}^{(j)}\left(t\right)\rtri|_{0, \left[t_{i-1}, t_{i}\right)}\leq \left[1+Eh^{2}\right]\left\|\overrightarrow{u}^{(j)}\left(t_{i-1}\right)\right\|+hD\ltri|F^{(j)}\left(t\right)\rtri|_{0},
\end{equation}
\begin{equation}\label{Chapt_3_2_eq_66}
    \left\|\overrightarrow{u}^{(j)}\left(t_{i}\right)\right\|\leq \left[e^{-\alpha h_{i}}+Eh_{i}^{2}\right]\left\|\overrightarrow{u}^{(j)}\left(t_{i-1}\right)\right\|+h_{i}D\ltri|F^{(j)}\left(t\right)\rtri|_{0},
\end{equation}
where
$$E=\cfrac{B+CP}{2},\quad D=1+\cfrac{CQ\overline{h}}{2}.$$
Let us make in \eqref{Chapt_3_2_eq_66} the following substitution
\begin{equation}\label{Chapt_3_2_eq_67}
    \left\|\overrightarrow{u}^{(j)}\left(t_{i}\right)\right\|\ltri|F^{(j)}\left(t\right)\rtri|_{0}^{-1}= y_{i},\quad \ltri|F^{(j)}\left(t\right)\rtri|_{0}>0,
\end{equation}
then the system of inequalities \eqref{Chapt_3_2_eq_66} turns to the following one
\begin{equation}\label{Chapt_3_2_eq_68}
    y_{i}\leq \left[e^{-\alpha h_{i}}+Eh_{i}^{2}\right]y_{i-1}+h_{i}D,\quad y_{0}=0,\; i\in \N.
\end{equation}
Now we are going to prove that for every sequence of sufficiently small constants $h_{i}>0,\; i\in \N$ a sequence  $\left\{y_{i}\right\}_{i=0}^{\infty}$ which satisfies  inequalities \eqref{Chapt_3_2_eq_68} is bounded. So far as for all $h_{i}\geq 0$ the following inequality holds true
$$e^{-\alpha h_{i}}+Eh_{i}^{2}\leq 1-\alpha h_{i}+h_{i}^{2}\bigg(E+\cfrac{\alpha^{2}}{2}\bigg)=1-h_{i}\bigg(\alpha -h_{i}\bigg(E+\cfrac{\alpha^{2}}{2}\bigg)\bigg),$$
the assumption
\begin{equation}\label{Chapt_3_2_eq_69}
    h_{i}\leq\cfrac{\alpha}{2E+\alpha^{2}},\quad i\in \N,
\end{equation}
yields the inequality
$$e^{-\alpha h_{i}}+Eh_{i}^{2}\leq 1-\alpha h_{i} + h_{i}^{2}\bigg(E+\cfrac{\alpha^{2}}{2}\bigg)\leq 1-\cfrac{h_{i}\alpha}{2}.$$
Taking into account the last estimate it is easy to verify that under conditions \eqref{Chapt_3_2_eq_69} the recursive sequence
\begin{equation}\label{Chapt_3_2_eq_70}
    Y_{i}=\bigg(1-h_{i}\cfrac{\alpha}{2}\bigg)Y_{i-1}+h_{i}D,\quad i\in \N,\; Y_{0}=0,
\end{equation}
is dominant for every sequence $\left\{y_{i}\right\}_{i=0}^{\infty}$ which satisfies inequalities \eqref{Chapt_3_2_eq_68}. If we make the substitution
$$h_{i}=\mathfrak{h}_{i}\cfrac{2}{\alpha},\; Y_{i}=\cfrac{2 D}{\alpha}\left(1-Z_{i}\right),\quad i\in \N,$$ in  \eqref{Chapt_3_2_eq_70} we would obtain the following recursive sequence
$$Z_{i}=\left(1-\mathfrak{h}_{i}\right)Z_{i-1},\quad i\in \N,\; Z_{0}=1.$$
Thus
$$Z_{i}=\prod\limits_{p=1}^{i}\left(1-\mathfrak{h}_{p}\right),\quad i\in \N, \; Z_{0}=1,$$
and the assumptions
\begin{equation}\label{Chapt_3_2_eq_71}
    \forall h_{i} \colon 0<h_{i}\leq \widetilde{h} =\min \bigg\{\cfrac{4}{\alpha}, \cfrac{\alpha}{2E+\alpha^{2}}, \overline{h}\left(\varepsilon_{1}\right)\bigg\},
\end{equation}
implies the estimate
\begin{equation}\label{Chapt_3_2_eq_72}
    0\leq y_{i}\leq Y_{i}\leq \cfrac{4}{\alpha} D.
\end{equation}
Using estimate \eqref{Chapt_3_2_eq_72} and taking into account \eqref{Chapt_3_2_eq_67} from inequalities \eqref{Chapt_3_2_eq_64}, \eqref{Chapt_3_2_eq_65} we can get
\begin{equation}\label{Chapt_3_2_eq_73}
    \ltri|\cfrac{d}{d t}\overrightarrow{u}^{(j)}\left(t\right)\rtri|_{0}\leq \bigg(\cfrac{4PD}{\alpha}+Q\bigg)\ltri|F^{(j)}\left(t\right)\rtri|_{0},
\end{equation}
\begin{equation}\label{Chapt_3_2_eq_74}
    \ltri|\overrightarrow{u}^{(j)}\left(t\right)\rtri|_{0}\leq \bigg(\left(1+E\widetilde{h}^{2}\right)\cfrac{4D}{\alpha}+\widetilde{h}D\bigg)\ltri|F^{(j)}\left(t\right)\rtri|_{0}.
\end{equation}
Let us put
\begin{equation*}\label{Chapt_3_2_eq_74'}
    \sigma = \max\bigg\{\cfrac{4PD}{\alpha}+Q,\;\left(1+E\widetilde{h}^{2}\right)\cfrac{4D}{\alpha}+\widetilde{h}D\bigg\},
\end{equation*}
  then from \eqref{Chapt_3_2_eq_73} and \eqref{Chapt_3_2_eq_74} we obtain
\begin{equation}\label{Chapt_3_2_eq_75}
    \ltri|\overrightarrow{u}^{(j)}\left(t\right)\rtri|\leq\sigma \ltri|F^{(j)}\left(t\right)\rtri|_{0}.
\end{equation}

In the following inequalities we use the assertions of corollary \ref{N_z lemi pro Adom1} and lemma \ref{L_pro Adom2} because it is easy to see that their assumptions are satisfied if we put
\begin{equation}\label{Chapt_3_2_eq_75'_pravki_general}
\widetilde{N}\left(u\right)=\sum\limits_{i=0}^{\infty}B_{i}u^{i},\; \forall u\in \R.
\end{equation}
 Indeed, using the Leibniz formula for the $n$-th derivative of the product (see, for example \cite[p. 268--272]{shvartz_analiz}) we can check that the function $\widetilde{N}\left(t\right)$ \eqref{Chapt_3_2_eq_75'_pravki_general} satisfies conditions \eqref{Chapt_3_2_eq_371}:
$$\begin{array}{c}
    \left.\cfrac{d^{k}B_{p}u^{p}}{du^{k}}\right|_{u=\left\|\overrightarrow{v}\right\|}\equiv\left.\left(B_{p}u^{p}\right)^{(k)}\right|_{u=\left\|\overrightarrow{v}\right\|}\geq \\[1.2em]
     \left.\geq\bigg(\sum\limits_{i_{1}+\ldots+i_{m}=p}\ltri|\mathbf{N}_{i_{1}\ldots i_{m}}\left(t\right)\rtri|_{0}u^{i_{1}+\ldots+ i_{m}}\bigg)^{(k)}\right|_{u=\left\|\overrightarrow{v}\right\|}= \\ [1.2em]
    = \bigg.\sum\limits_{k_{1}+\ldots+k_{m}=k}\cfrac{k!}{k_{1}!\ldots k_{m}!}\sum\limits_{i_{1}+\ldots+i_{m}=p}\ltri|\mathbf{N}_{i_{1}\ldots i_{m}}\left(t\right)\rtri|_{0}\left(u^{i_{1}}\right)^{(k_{1})}\ldots \left(u^{i_{m}}\right)^{(k_{m})}\bigg|_{u=\left\|\overrightarrow{v}\right\|}= \\ [1.2em]
    =\bigg.\sum\limits_{k_{1}+\ldots+k_{m}=k}\cfrac{k!}{k_{1}!\ldots k_{m}!}\sum\limits_{i_{1}+\ldots+i_{m}=p}\ltri|\mathbf{N}_{i_{1}\ldots i_{m}}\left(t\right)\rtri|_{0} \cfrac{\partial^{k} \left(u_{1}^{i_{1}}u_{2}^{i_{2}}\ldots u_{m}^{i_{m}}\right)}{\partial u_{1}^{k_{1}}\partial u_{2}^{k_{2}}\ldots \partial u_{m}^{k_{m}}}\bigg|_{                                                                   u_{j}=\left\|\overrightarrow{v}\right\|}\geq \\ [1.2em]
    \geq\bigg. \bigg.\sum\limits_{k_{1}+\ldots+k_{m}=k}\cfrac{k!}{k_{1}!\ldots k_{m}!}\left\|\sum\limits_{i_{1}+\ldots+i_{m}=p}\cfrac{\partial^{k}\; \mathbf{N}_{i_{1}\ldots i_{m}}\left(t\right)u_{1}^{i_{1}}u_{2}^{i_{2}}\ldots u_{m}^{i_{m}}}{\partial u_{1}^{k_{1}}\partial u_{2}^{k_{2}}\ldots \partial u_{m}^{k_{m}}}\right\|\bigg|_{u_{j}=v_{j}},\\[1.2em]
    \forall\; \overrightarrow{v}=\left[v_{1}, v_{2}, \ldots, v_{m}\right]^{T}\in V_{m}\left(\R\right).
    \end{array}
$$

From \eqref{Chapt_3_2_eq_75} we have
\begin{equation}\label{Chapt_3_2_eq_76}
    \begin{array}{c}
      \ltri|\overrightarrow{u}^{(j)}\left(t\right)\rtri|\leq \sigma \Bigg\{\sum\limits_{p=1}^{j-1}A_{j-p}\left(\widetilde{N}\left(\cdot\right); \left[\ltri|\overrightarrow{u}^{(i)}\left(t\right)\rtri|\right]_{i=0}^{j-p}\right)\ltri|\overrightarrow{u}^{(p)}\left(t\right)\rtri|+\Bigg. \\[1.2em]
      +h\sum\limits_{p=0}^{j-1}A_{j-1-p}\left(\widetilde{N}^{\prime}\left(\cdot\right)\times \left(\cdot\right); \left[\ltri|\overrightarrow{u}^{(i)}\left(t\right)\rtri|\right]_{i=0}^{j-1-p}\right)\ltri|\overrightarrow{u}^{(p)}\left(t\right)\rtri|+ \\[1.2em]
      +\cfrac{\ltri|\overrightarrow{u}^{(0)}\left(t\right)\rtri|}{j!}\cfrac{d^{j}}{d \tau^{j}}\Bigg[\widetilde{N}\left(\sum\limits_{i=0}^{\infty}\tau^{i}\ltri|\overrightarrow{u}^{(i)}\left(t\right)\rtri|\right)-\Bigg.\\[1.2em]
      -\Bigg.\bigg.\Bigg.\bigg.\cfrac{d \widetilde{N}\left(u\right)}{d u}\bigg|_{u=\ltri|\overrightarrow{u}^{(0)}\left(t\right)\rtri|}\tau^{j}\ltri|\overrightarrow{u}^{(j)}\left(t\right)\rtri|\Bigg]\bigg|_{\tau=0}\Bigg\},\quad j\in \N.
    \end{array}
\end{equation}
Let us put
\begin{equation}\label{Chapt_3_2_eq_77}
    h^{-j}\ltri|\overrightarrow{u}^{(j)}\left(t\right)\rtri|=v_{j},\quad j\in \N.
\end{equation}
in \eqref{Chapt_3_2_eq_76}.
Then, changing $v_{j}$ by  $V_{j}$ and the inequality sign to the equality one we obtain the following system of equations with respect to $V_{j}$
\begin{equation*}\label{Chapt_3_2_eq_78}
    \begin{array}{c}
      V_{j}=\sigma \Bigg\{\sum\limits_{p=1}^{j-1}A_{j-p}\left(\widetilde{N}\left(\cdot\right); \left[V_{i}\right]_{i=0}^{j-p}\right)V_{p}+\Bigg.\sum\limits_{p=0}^{j-1}A_{j-1-p}\left(\widetilde{N}^{\prime}\left(\cdot\right)\times \left(\cdot\right); \left[V_{i}\right]_{i=0}^{j-1-p}\right)V_{p}+ \\[1.2em]
      \Bigg.\bigg.+\cfrac{V_{0}}{j!}\cfrac{d^{j}}{d \tau^{j}}\bigg(\widetilde{N}\bigg(\sum\limits_{p=0}^{\infty}\tau^{p}V_{p}\bigg)\bigg)\bigg|_{\tau=0}-V_{j}V_{0}\widetilde{N}^{\prime}\left(V_{0}\right)\Bigg\},\; j\in \N,
    \end{array}
\end{equation*}
where $V_{0}=\mu,$ or
\begin{equation}\label{Chapt_3_2_eq_79}
    \begin{array}{c}
      V_{j}=\cfrac{\sigma}{1+\sigma V_{0}\widetilde{N}^{\prime}\left(V_{0}\right)}\Bigg\{\sum\limits_{p=0}^{j-1}A_{j-p}\left(\widetilde{N}\left(\cdot\right); \left[V_{i}\right]_{i=0}^{j-p}\right)V_{p}+\Bigg. \\[1.2em]
      +\sum\limits_{p=0}^{j-1}A_{j-1-p}\left(\widetilde{N}^{\prime}\left(\cdot\right)\times \left(\cdot\right); \left[V_{i}\right]_{i=0}^{j-1-p} \right)V_{p}\Bigg\},\;j\in \N.
    \end{array}
\end{equation}
The sequence $\left\{V_{i}\right\}_{i=1}^{\infty},$ which satisfies the system of recursive equations \eqref{Chapt_3_2_eq_79}, is dominant for the sequence $\left\{v_{i}\right\}_{i=1}^{\infty}.$ Using the generating functions method, we are going to prove that for the $h$ sufficiently small the assumptions \eqref{Chapt_3_2_eq_55}, \eqref{Chapt_3_2_eq_56} hold and the power series $\sum\limits_{i=0}^{\infty}z^{i}V_{i}$ possess a nonzero convergent radius.

From \eqref{Chapt_3_2_eq_79} we have
\begin{equation}\label{Chapt_3_2_eq_81}
    \begin{array}{c}
      g\left(z\right)-V_{0}=\cfrac{\sigma}{1+\sigma V_{0}\widetilde{N}^{\prime}\left(V_{0}\right)}\left\{g\left(z\right)\left(\widetilde{N}\left(g\left(z\right)\right)-\widetilde{N}\left(V_{0}\right)\right)\right.+ \\[1.2em]
      \left.+zg^{2}\left(z\right)\widetilde{N}^{\prime}\left(g\left(z\right)\right)\right\}, \end{array}
\end{equation}
where
\begin{equation}\label{Chapt_3_2_eq_80}
    g\left(z\right)=\sum\limits_{j=0}^{\infty}z^{j}V_{j}.
\end{equation}

To prove that the function $g\left(z\right),$ which satisfies equality \eqref{Chapt_3_2_eq_81}, exists or, in other words, has nonempty domain which includes some neighbourhood of the point $z=0,$ we have to express $z$ from \eqref{Chapt_3_2_eq_81}  by means of $g \colon$
\begin{equation}\label{Chapt_3_2_eq_82}
    \begin{array}{c}
      z\left(g\right)=\cfrac{1}{g^{2}\widetilde{N}^{\prime}\left(g\right)}\left\{\Sigma\left(g-V_{0}\right)+\left(\widetilde{N}\left(V_{0}\right)-\widetilde{N}\left(g\right)\right)g\right\}, \\[1.2em]
       V_{0}\leq g,\; \Sigma=\cfrac{1}{\sigma}+V_{0}\widetilde{N}^{\prime}\left(V_{0}\right).
    \end{array}
\end{equation}
It is easy to see that the function $z\left(g\right)$ \eqref{Chapt_3_2_eq_82} is defined and continuously differentiable in some neighbourhood of the point $g=V_{0}.$ To prove that there exists an inverse function  $g=g\left(z\right)$ defined in some neighbourhood of the point $z=0,$ it is sufficient to check that $z^{\prime}\left(V_{0}\right)>0.$ The last inequality follows from the expression for $z\left(g\right)$ \eqref{Chapt_3_2_eq_82} $\colon$
\begin{equation}\label{Chapt_3_2_eq_83}
    \begin{array}{c}
      z^{\prime}\left(V_{0}\right)=\lim\limits_{g\rightarrow V_{0}}\cfrac{z\left(g\right)-z\left(V_{0}\right)}{g-V_{0}}= \\[1.2em]
      =\lim\limits_{g\rightarrow V_{0}}\cfrac{1}{g^{2}\widetilde{N}^{\prime}\left(g\right)}\Bigg(\Sigma-g\cfrac{\widetilde{N}\left(g\right)-\widetilde{N}\left(V_{0}\right)}{g-V_{0}}\Bigg)=\cfrac{1}{\sigma V_{0}^{2}\widetilde{N}^{\prime}\left(V_{0}\right)}>0.
    \end{array}
\end{equation}

 It follows from inequality \eqref{Chapt_3_2_eq_83} (see \cite[p. 87]{An_theory}) that there exists an inverse function $f\left(z\right)$ which is holomorphic in  some open neighbourhood of the point $z=0:$ $\mathcal{C}\subseteq \mathbb{C}.$ Assume that $R_{1}>0$ is the largest constant for which the open ball $B_{R_{1}}=\{z\in \mathbb{C} \colon |z|<R_{1}\}$ satisfies inclusion $B_{R_{1}}\subseteq \mathcal{C}.$ Now we are going to prove that series \eqref{Chapt_3_2_eq_80} converges at the point $z=R_{1}.$ Suppose that this is not true and series \eqref{Chapt_3_2_eq_80} diverges at the point $z=R_{1}.$ Then we have $$\lim\limits_{z\to R_{1}-0}f(z)=+\infty.$$ But, on the other hand, for all $z$ from $(-R_{1}, R_{1})$ equality \eqref{Chapt_3_2_eq_81} holds. And we immediately get the following contradiction
 \begin{eqnarray}
   1 = f(z)-\overline{v}_{0}+ \\
     + \frac{(1+\overline{v}_0)z[\|q'\|_{0,1}f(z)+\overline{N}(f(z))+\overline{N}^{\prime}\left(\overline{v}_{0}\right)\overline{v}_{0}-\overline{N}\left(\overline{v}_{0}\right)]}{f(z)-\overline{v}_{0}}\to + \infty,
 \end{eqnarray}
 as $z \to R_{1}-0.$
   This contradiction implies the inequality $f(R_{1})<+\infty$ and
$$R_{1}^{j}V_{j}\leq\cfrac{C}{\left(j+1\right)^{1+\varepsilon}}$$ for some positive constants $\varepsilon$ and $C,$ which do not depend on $h.$

If we return to notations \eqref{Chapt_3_2_eq_77} we would obtain an estimate
\begin{equation}\label{Chapt_3_2_eq_84}
    \ltri|\overrightarrow{u}^{(j)}\left(t\right)\rtri|\leq \cfrac{C}{\left(j+1\right)^{1+\varepsilon}}\bigg(\cfrac{h}{R_{1}}\bigg)^{j}, \; j\in \N\bigcup\left\{0\right\},
\end{equation}
and the following sufficient condition for series \eqref{Chapt_3_2_eq_80} be convergent
\begin{equation*}\label{Chapt_3_2_eq_84'}
  \cfrac{h}{R_{1}}\leq 1.
\end{equation*}

So, we have proved that for the sufficiently small maximum step $h>0$ of the grid $\widehat{\omega}$ \eqref{Chapt_3_1_eq_4},  namely,
\begin{equation}\label{Chapt_3_2_eq_84'_general_pravki_1}
h\leq R=\min\left\{\widetilde{h}, R_{1}\right\}
\end{equation}
 (see \eqref{Chapt_3_2_eq_71}) assumptions \eqref{Chapt_3_2_eq_55}, \eqref{Chapt_3_2_eq_56} hold true. Hence, in accordance with definition \ref{O_FD_zbizhn}, the FD-method for the Cauchy problem \eqref{Chapt_3_1_eq_1} converges. Now we have to prove that it converges to the exact solution of problem \eqref{Chapt_3_1_eq_1} by the same definition.

 Let us assume that the grig $\widehat{\omega}$ \eqref{Chapt_3_1_eq_4} satisfies condition \eqref{Chapt_3_2_eq_84'_general_pravki_1}. We define (see notations on page \pageref{P_funk_vect_norm})
\begin{equation}\label{Chapt_3_2_eq_86}
    \overrightarrow{v}\left(t\right)=\sum\limits_{j=0}^{\infty}\overrightarrow{u}^{(j)}\left(t\right) \in V_{m}\left(\mathbb{Q}^{1}_{\omega}\left[t_{0}, +\infty\right)\right).
\end{equation}
It easy to see that
 $$\cfrac{d}{d t}\overrightarrow{v}\left(t\right)=\sum\limits_{j=0}^{\infty}\cfrac{d}{d t}\overrightarrow{u}^{(j)}\left(t\right),\; \forall t\in \left(t_{i-1}, t_{i}\right),\; \forall i\in\N.$$
Let us prove that vector-function \eqref{Chapt_3_2_eq_86} is the solution of the Cauchy problem \eqref{Chapt_3_1_eq_1}. For this reason we need to sum the equation of the base problem \eqref{Chapt_3_1_eq_5} with equations \eqref{Chapt_3_1_eq_6} $\forall j\in \N.$ It results in the following equalities
\begin{equation}\label{Chapt_3_2_eq_85}
\begin{array}{c}
  \cfrac{d}{d t}\overrightarrow{v}\left(t\right)-\sum\limits_{j=0}^{\infty}\sum\limits_{p=0}^{j}A_{j-p}\left(\mathbf{N}\left(t, \cdot\right); \left[\overrightarrow{u}^{(i)}\left(t\right)\right]_{i=0}^{j-p}\right)\overrightarrow{u}^{(p)}\left(t\right)=\overrightarrow{\phi}\left(t\right), \\ [1.2em]
  \cfrac{d}{d t}\overrightarrow{v}\left(t_{i}+0\right)-\sum\limits_{j=0}^{\infty}\sum\limits_{p=0}^{j}A_{j-p}\left(\mathbf{N}\left(t_{i}, \cdot\right); \left[\overrightarrow{u}^{(i)}\left(t_{i}\right)\right]_{i=0}^{j-p}\right)\overrightarrow{u}^{(p)}\left(t_{i}\right)=\overrightarrow{\phi}\left(t_{i}\right), \end{array}
\end{equation}
$\forall i\in\N\bigcup\left\{0\right\},\;\forall t\in \bigcup\limits_{i=1}^{\infty}\left(t_{i-1}, t_{i}\right).$

Using the theorem about the substitution of the power series into the power series (see, for example, \cite[p. 485]{Fihtenholts}), it is easy to obtain the following equality
$$\mathbf{N}\bigg(t, \sum\limits_{i=0}^{\infty}\tau^{i}\overrightarrow{u}^{(i)}\left(t\right)\bigg)\sum\limits_{i=0}^{\infty}\tau^{i}\overrightarrow{u}^{(i)}\left(t\right)=$$
$$=\sum\limits_{i=0}^{\infty}\tau^{j}\sum\limits_{p=0}^{j}A_{j-p}\left(\mathbf{N}\left(t, \cdot\right); \left[\overrightarrow{u}^{(i)}\left(t\right)\right]_{i=0}^{j-p}\right)\overrightarrow{u}^{(p)}\left(t\right),$$
$\forall t\in \left[0, 1\right],$ $\forall t\in \left[t_{0}, +\infty\right).$ Taking it into consideration, we can rewrite equalities \eqref{Chapt_3_2_eq_85} in the following form
\begin{equation}\label{Chapt_3_2_eq_87}
\begin{array}{c}
  \cfrac{d}{d t}\overrightarrow{v}\left(t\right)-\mathbf{N}\left(t, \overrightarrow{v}\left(t\right)\right)\overrightarrow{v}\left(t\right)=\overrightarrow{\phi}\left(t\right), t\in \bigcup\limits_{i=1}^{\infty}\left(t_{i}, t_{i-1}\right), \\ [1.2em]
   \cfrac{d}{d t}\overrightarrow{v}\left(t_{i}+0\right)-\mathbf{N}\left(t_{i}, \overrightarrow{v}\left(t_{i}\right)\right)\overrightarrow{v}\left(t_{i}\right)=\overrightarrow{\phi}\left(t_{i}\right),\; i\in\N\bigcup\left\{0\right\}.
\end{array}
\end{equation}
The uniqueness of the solution $\overrightarrow{u}\left(t\right)$ of the Cauchy problem \eqref{Chapt_3_1_eq_1} together with equalities \eqref{Chapt_3_2_eq_87} imply that as far as $\overrightarrow{v}\left(t_{i-1}\right)=\overrightarrow{u}\left(t_{i-1}\right)$ the identity $\overrightarrow{v}\left(t\right)\equiv\overrightarrow{u}\left(t\right),$ $\forall t\in\left[t_{i-1}, t_{i}\right],$ $i\in\N$ holds true. Thus, the evident equality $\overrightarrow{v}\left(t_{0}\right)=\overrightarrow{u}\left(t_{0}\right)=\overrightarrow{u}_{0}$ implies the identity $\overrightarrow{v}\left(t\right)\equiv\overrightarrow{u}\left(t\right),$ $\forall t\in\left[t_{0}, +\infty\right).$

Combining \eqref{Chapt_3_2_eq_84} with \eqref{Chapt_3_2_eq_84'_general_pravki_1} it is easy to obtain estimates \eqref{Chapt_3_2_eq_18}, \eqref{Chapt_3_2_eq_19}.
 The theorem is proved. $\blacksquare$

{\bf Numerical example.}
Let us consider the following Cauchy problem
\begin{equation}\label{Chapt_3_5_eq_1}
\left\{
                    \begin{array}{c}
                      \cfrac{d}{d t}u_{1}\left(t\right)+\left(u_{1}^{2}\left(t\right)+u_{2}^{2}\left(t\right)+1\right)u_{1}\left(t\right)-u_{1}\left(t\right)u_{2}\left(t\right)=\\ =2\sin\left(t\right)-\cfrac{1}{2}\sin\left(2t\right)+\cos\left(t\right), \\ [1.2em]
                      \cfrac{d}{d t}u_{2}\left(t\right)+\left(u_{1}^{2}\left(t\right)+u_{2}^{2}\left(t\right)+1\right)u_{2}\left(t\right)-u_{1}\left(t\right)u_{2}\left(t\right)=\\
                      =2\cos\left(t\right)-\cfrac{1}{2}\sin\left(2t\right)-\sin\left(t\right). \\
                    \end{array}
                  \right.
\end{equation}
\begin{equation}\label{Chapt_3_5_eq_2}
    u_{1}\left(0\right)=0,\quad u_{2}\left(0\right)=1,\quad t\in\left[0, +\infty\right).
\end{equation}
Using the notations of system \eqref{Chapt_3_1_eq_1} we would have ($\overrightarrow{u}=\left[u_{1}, u_{2}\right]^{T}$) :
\begin{equation}\label{Chapt_3_5_eq_3}
\begin{array}{c}
  \mathbf{N}\left(\overrightarrow{u}\right)=\left[
                                                \begin{array}{cc}
                                                  u_{2}-u_{1}^{2}-u_{2}^{2}-1 & 0 \\
                                                  0 & u_{1}-u_{1}^{2}-u_{2}^{2}-1 \\
                                                \end{array}
                                              \right]= \\ [1.2em]
  =-\left[
      \begin{array}{cc}
        1 & 0 \\
        0 & 1 \\
      \end{array}
    \right]+u_{1}\left[
      \begin{array}{cc}
        0 & 0 \\
        0 & 1 \\
      \end{array}
    \right]+u_{2}\left[
      \begin{array}{cc}
        1 & 0 \\
        0 & 0 \\
      \end{array}
    \right]+u_{1}^{2}\left[
      \begin{array}{cc}
        1 & 0 \\
        0 & 1 \\
      \end{array}
    \right]+u_{2}^{2}\left[
      \begin{array}{cc}
        1 & 0 \\
        0 & 1 \\
      \end{array}
    \right],\\[1.2em]
\widetilde{N}\left(u\right)=1+2u+2u^{2},
\end{array}
\end{equation}
\begin{equation}\label{Chapt_3_5_eq_4}
    \overrightarrow{\phi}\left(t\right)=\left[
                                          \begin{array}{c}
                                            2\sin\left(t\right)-\cfrac{1}{2}\sin\left(2t\right)+\cos\left(t\right) \\
                                            2\cos\left(t\right)-\cfrac{1}{2}\sin\left(2t\right)-\sin\left(t\right) \\
                                          \end{array}
                                        \right].
\end{equation}
It is easy to verify that the Jacobian matrix of the vector-function
\begin{equation}\label{Chapt_3_5_eq_5}
    \mathfrak{N}\left(\overrightarrow{u}\right)=\mathbf{N}\left(\overrightarrow{u}\right)\overrightarrow{u}=\left[u_{1}\left(u_{2}-u_{1}^{2}-u_{2}^{2}-1\right), u_{2}\left(u_{1}-u_{1}^{2}-u_{2}^{2}-1\right) \right]^{T}
\end{equation}
has the form
\begin{equation*}\label{Chapt_3_5_eq_6}
    J\left(\overrightarrow{u}\right)=\left[
                                       \begin{array}{cc}
                                         u_{2}-u_{2}^{2}-3u_{1}^{2}-1 & u_{1}-2u_{1}u_{2} \\
                                         u_{2}-2u_{1}u_{2} & u_{1}-u_{1}^{2}-3u_{2}^{2}-1 \\
                                       \end{array}
                                     \right].
\end{equation*}
Applying the Sylvester’s criterion to the symmetric matrix $J_{s}\left(\overrightarrow{u}\right)=\cfrac{1}{2}\left(J\left(\overrightarrow{u}\right)+J^{T}\left(\overrightarrow{u}\right)\right),$ we obtain that the matrix $J_{s}\left(\overrightarrow{u}\right)$ is negative defined  $\forall \overrightarrow{u}\in V_{2}\left(\R\right)$. If we find the extremum values of the function $\cfrac{Tr\left(J_{s}\left(\overrightarrow{u}\right)\right)}{\det\left(J_{s}\left(\overrightarrow{u}\right)\right)},$ as the function of two variables $u_{1}, u_{2}$, we would get the inequality
\begin{equation}\label{Chapt_3_5_eq_7}
    0>\cfrac{Tr\left(J_{s}\left(\overrightarrow{u}\right)\right)}{\det\left(J_{s}\left(\overrightarrow{u}\right)\right)}=\cfrac{\lambda_{1}\left(\overrightarrow{u}\right)+\lambda_{2}\left(\overrightarrow{u}\right)}{\lambda_{1}\left(\overrightarrow{u}\right)\lambda_{2}\left(\overrightarrow{u}\right)}>-2.15,\quad \forall \overrightarrow{u}\in V_{2}\left(\R\right),
\end{equation}
 where $\lambda_{1}\left(\overrightarrow{u}\right), \lambda_{2}\left(\overrightarrow{u}\right)<0$ are the eigenvalues of the matrix $J_{s}\left(\overrightarrow{u}\right).$
From inequality \eqref{Chapt_3_5_eq_7} it follows that
\begin{equation}\label{Chapt_3_5_eq_8}
    \min\left\{\left|\lambda_{1}\left(\overrightarrow{u}\right)\right|, \left|\lambda_{2}\left(\overrightarrow{u}\right)\right|\right\}\geq \alpha=\left(2.15\right)^{-1}\approx 0.47.
\end{equation}

Similarly, if we find the extremum values of the function $\left\|\phi\left(t\right)\right\|,$ we would come to the estimate
\begin{equation}\label{Chapt_3_5_eq_9}
    \ltri|\phi\left(t\right)\rtri|_{0}<\kappa=2.9.
\end{equation}

So, as it was showed above, the Cauchy problem \eqref{Chapt_3_5_eq_1}, \eqref{Chapt_3_5_eq_2} satisfies the conditions of theorem \ref{T_first_FD} with the constants $\alpha$ \eqref{Chapt_3_5_eq_8} and $\kappa$ \eqref{Chapt_3_5_eq_9}. It is easy to see that the exact solution  $\overrightarrow{u}^{\ast}\left(t\right)$ of the problem is
\begin{equation}\label{Chapt_3_5_eq_10}
    \overrightarrow{u}^{\ast}\left(t\right)=\left[u^{\ast}_{1}\left(t\right), u^{\ast}_{2}\left(t\right)\right],\quad u^{\ast}_{1}\left(t\right)=\sin\left(t\right), \quad u_{2}^{\ast}\left(t\right)=\cos\left(t\right).
\end{equation}

As a first step we try to use the ADM to find the approximate solution of problem \eqref{Chapt_3_5_eq_1}, \eqref{Chapt_3_5_eq_2} (see, for example, \cite{Adomian_1994}).

\begin{figure}[htbp]
\begin{minipage}[h]{1\linewidth}
\begin{minipage}[h]{0.48\linewidth}
\center{\rotatebox{-0}{\includegraphics[height=0.9\linewidth,
width=1.3\linewidth]{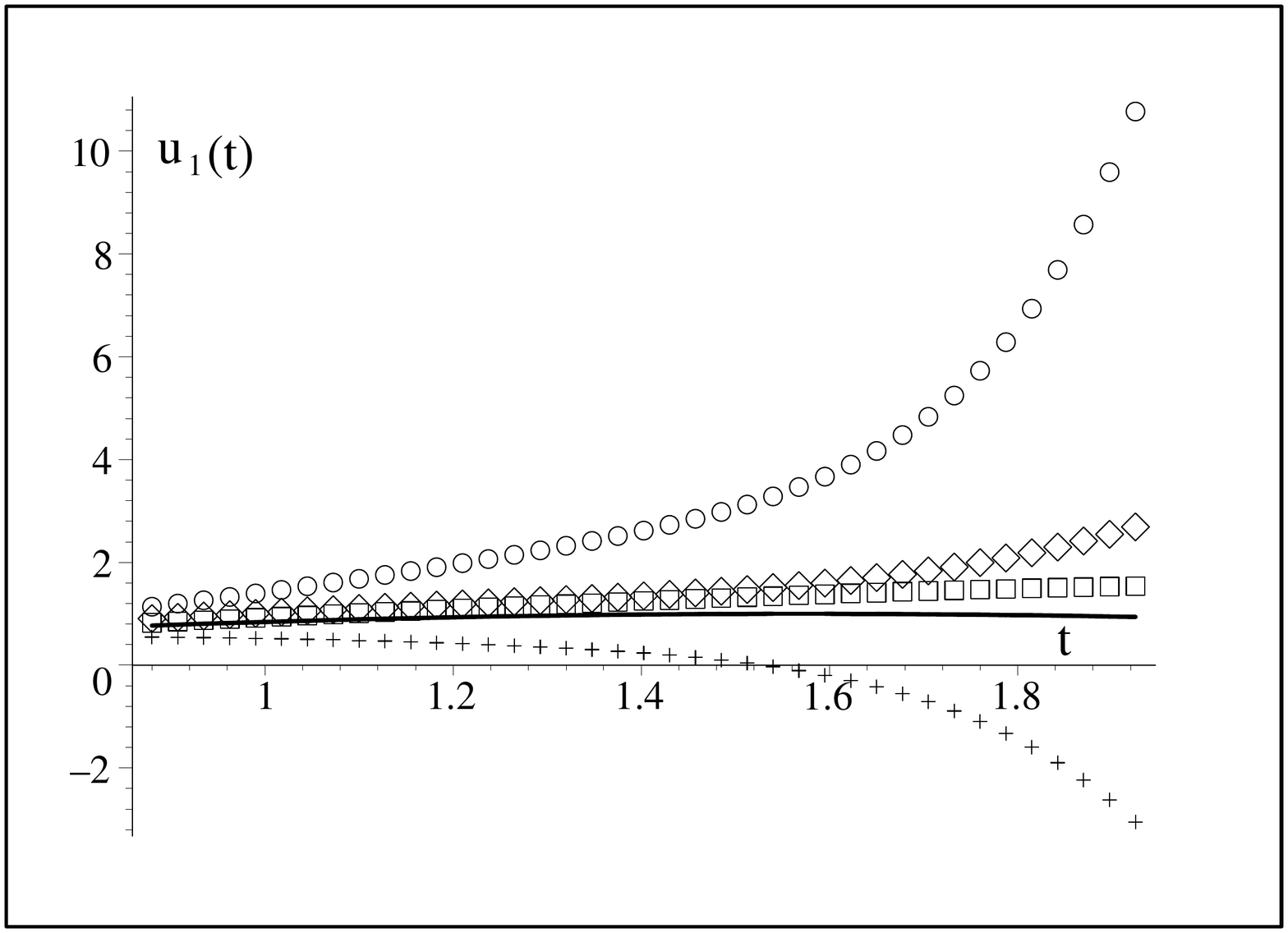}} \\ a) }
\end{minipage}
\hfill
\begin{minipage}[h]{0.48\linewidth}
\center{\rotatebox{-0}{\includegraphics[height=0.9\linewidth,
width=1.3\linewidth]{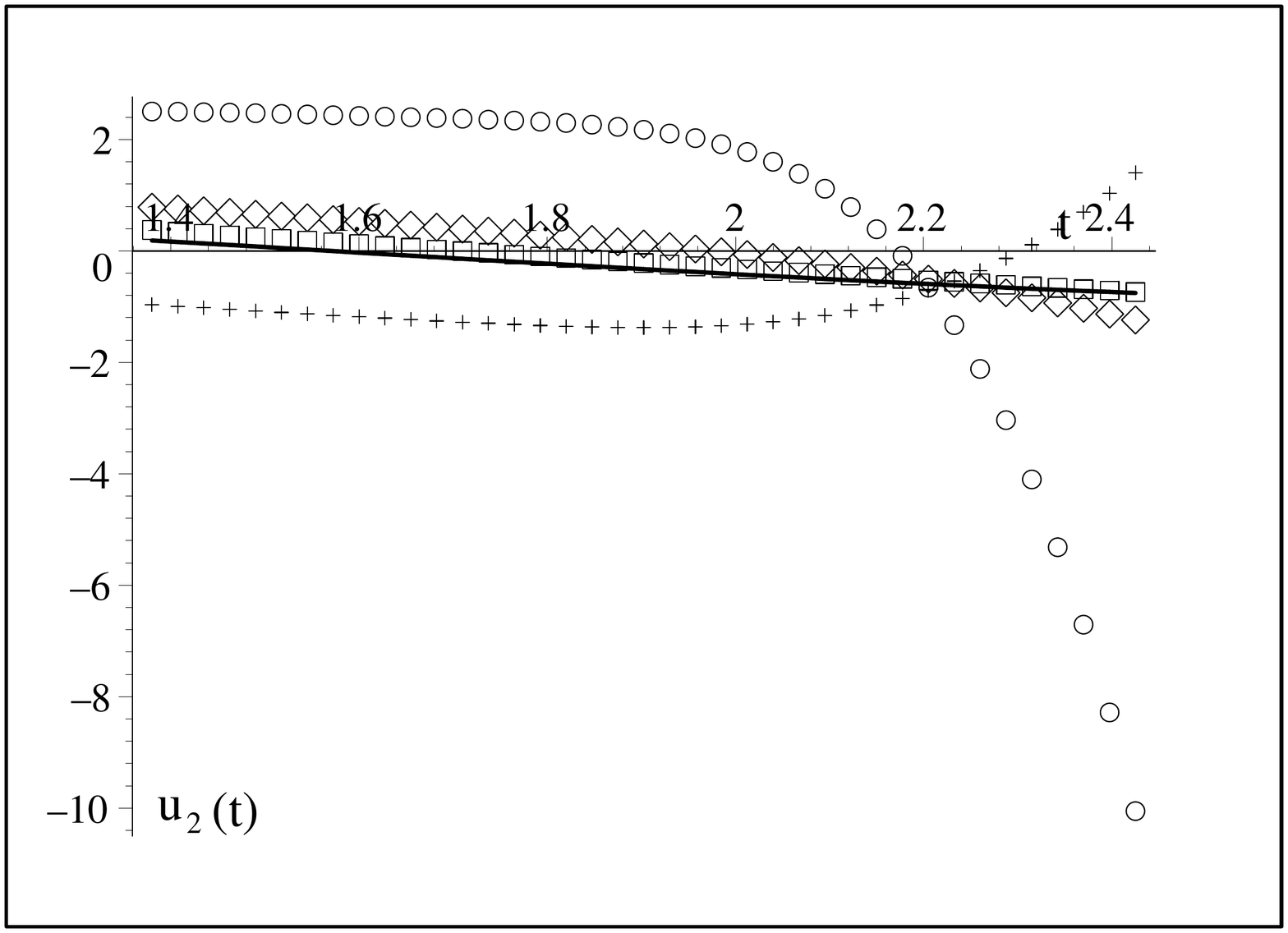}} \\ b)}
\end{minipage}
\caption{The Cauchy problem \eqref{Chapt_3_5_eq_1}, \eqref{Chapt_3_5_eq_2}, ADM. $\Box\colon \overset{0}{u}_{A i}\left(t\right),$ $\diamondsuit\colon \overset{2}{u}_{A i}\left(t\right),$ $+\colon \overset{3}{u}_{A i}\left(t\right),$ $\circ\colon \overset{4}{u}_{A i}\left(t\right),$ the solid line represents the exact solution \eqref{Chapt_3_5_eq_10},  a) $i=1$, b) $i=2.$  }\label{image2}
\end{minipage}
\end{figure}

\begin{figure}[htbp]
\begin{minipage}[h]{1\linewidth}
\begin{minipage}[h]{0.48\linewidth}
\center{\rotatebox{-0}{\includegraphics[height=0.9\linewidth,
width=1.3\linewidth]{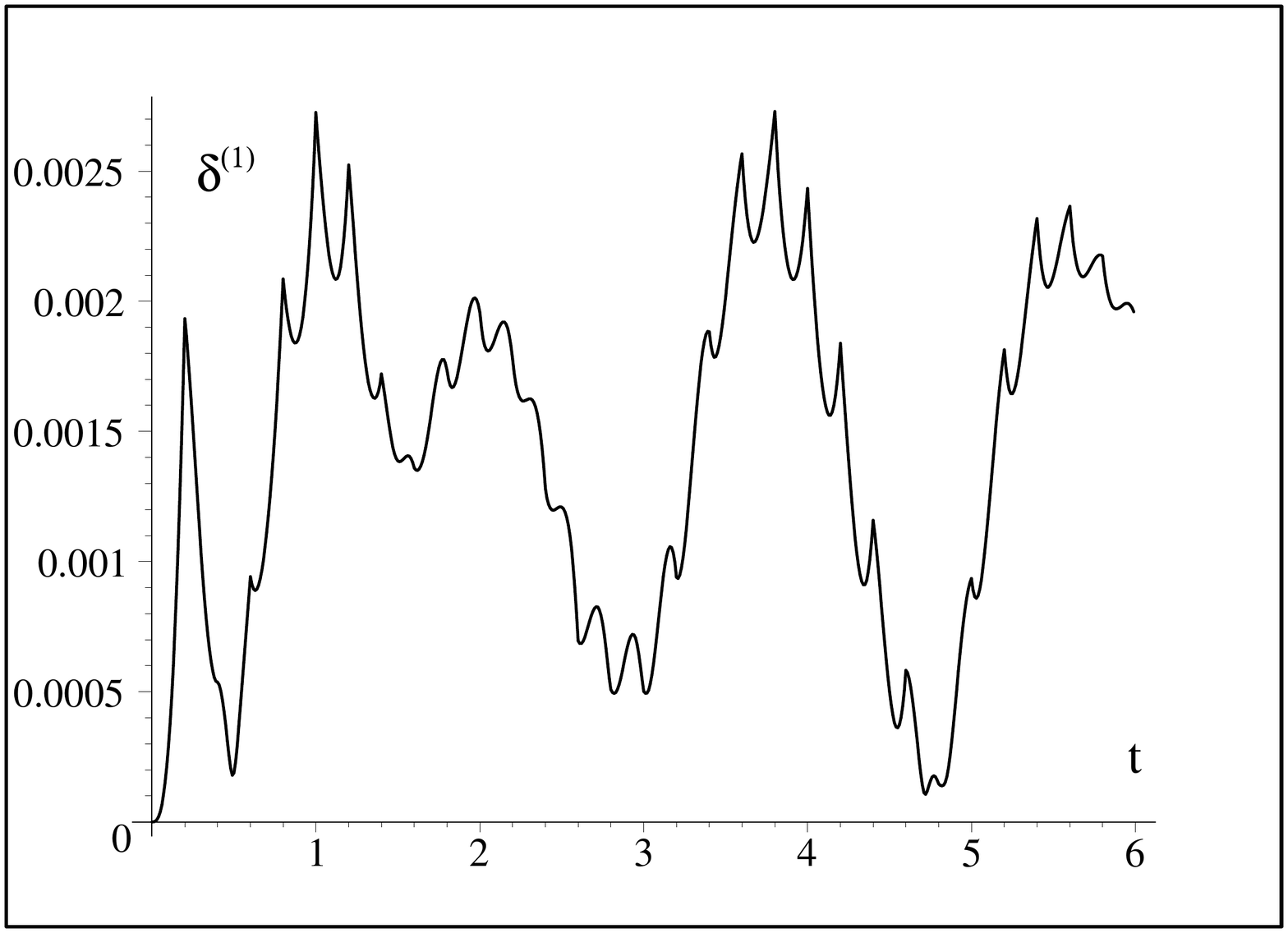}} \\ a) }
\end{minipage}
\hfill
\begin{minipage}[h]{0.48\linewidth}
\center{\rotatebox{-0}{\includegraphics[height=0.9\linewidth,
width=1.3\linewidth]{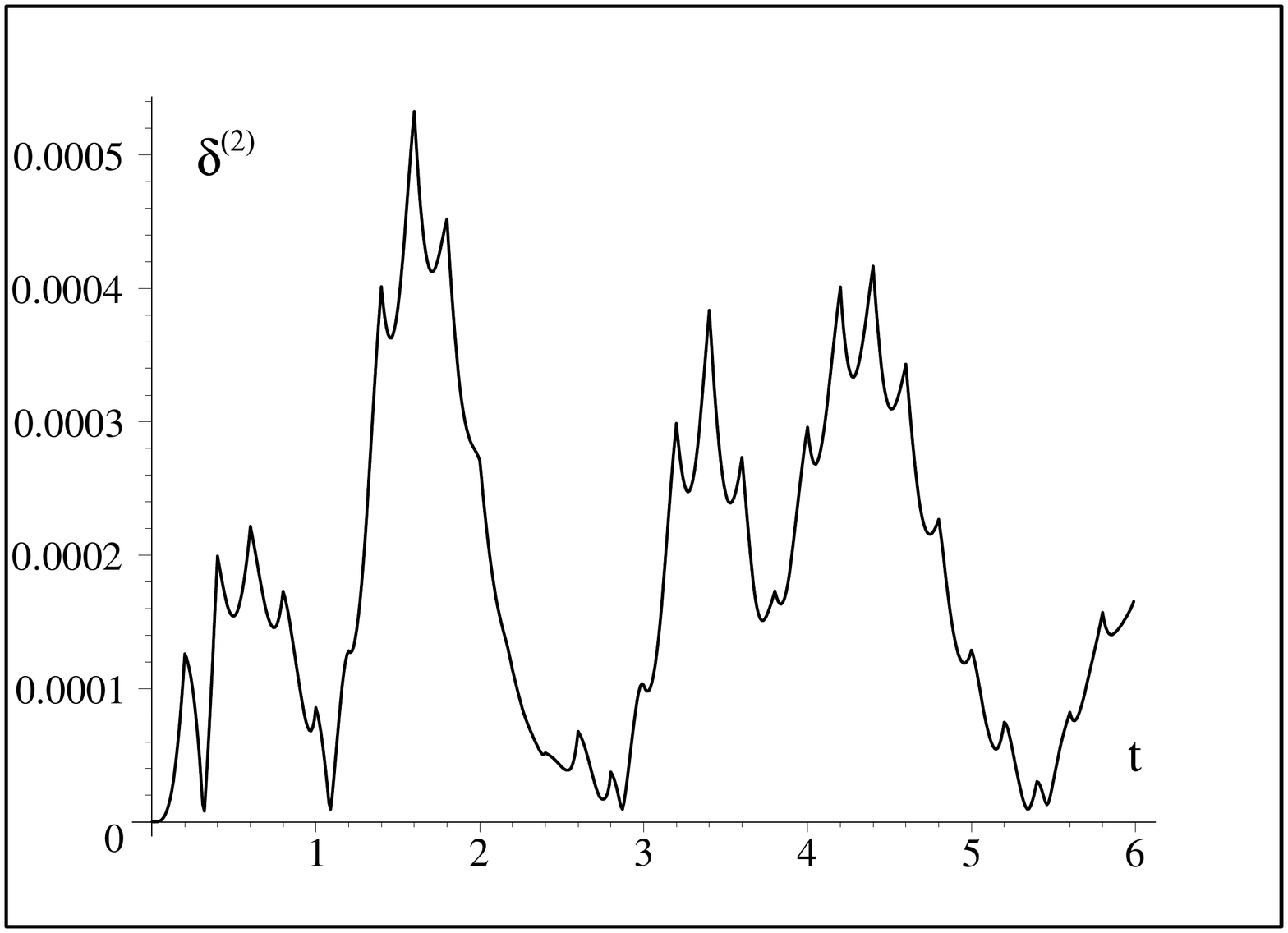}} \\ b)}
\end{minipage}
\begin{minipage}[h]{0.48\linewidth}
\center{\rotatebox{-0}{\includegraphics[height=0.9\linewidth,
width=1.3\linewidth]{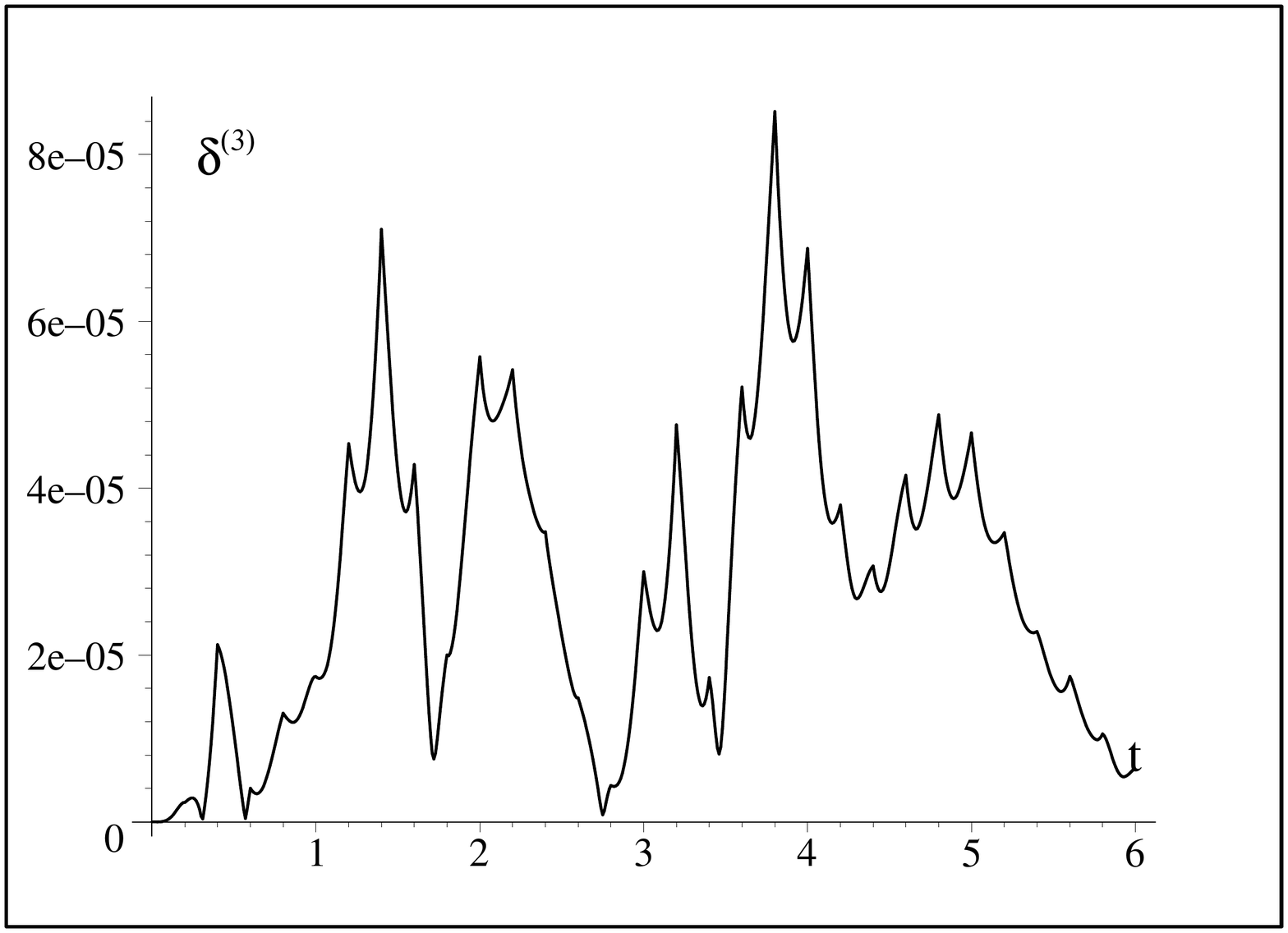}} \\ c) }
\end{minipage}
\hfill
\begin{minipage}[h]{0.48\linewidth}
\center{\rotatebox{-0}{\includegraphics[height=0.9\linewidth,
width=1.3\linewidth]{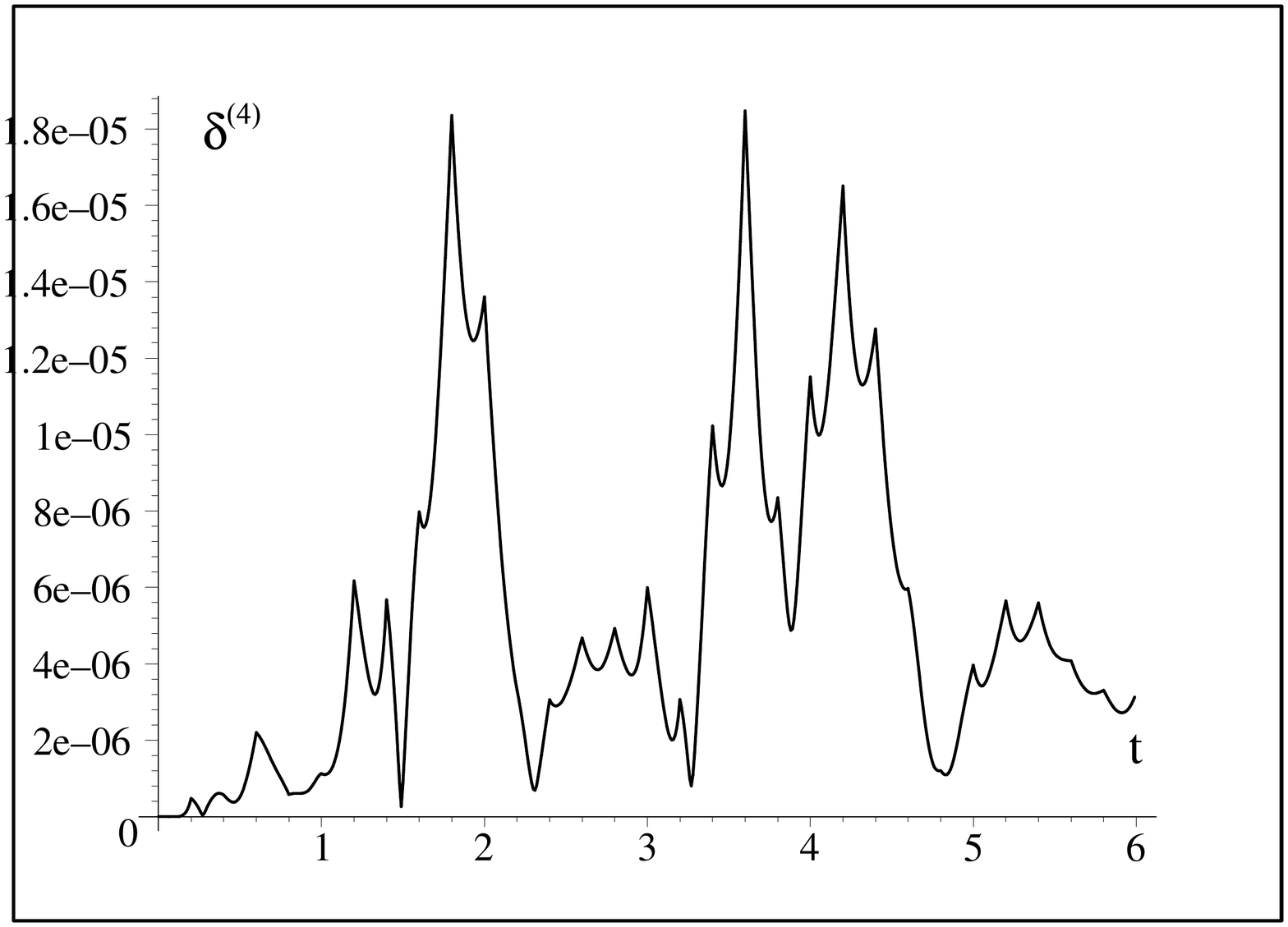}} \\ d)}
\end{minipage}
\caption{The Cauchy problem \eqref{Chapt_3_5_eq_1}, \eqref{Chapt_3_5_eq_2}, the FD-method. a) $\delta^{(1)}\left(t\right),$ b) $\delta^{(2)}\left(t\right),$\newline c) $\delta^{(3)}\left(t\right),$ d) $\delta^{(4)}\left(t\right).$ }\label{image3}
\end{minipage}
\end{figure}

For this reason we  need to express the vector-function $\mathfrak{N}\left(\overrightarrow{u}\right)$ \eqref{Chapt_3_5_eq_5} as a sum of a linear component and an essentially nonlinear one: \mbox{$\mathfrak{N}\left(\overrightarrow{u}\right)=\mathfrak{N}_{1}\left(\overrightarrow{u}\right)-E^{(2)}\overrightarrow{u}$.} We looking for the $p$-th approximation of the exact solution of the Cauchy problem \eqref{Chapt_3_5_eq_1}, \eqref{Chapt_3_5_eq_2} in the form of the partial sum $\overset{p}{\overrightarrow{u}}_{A}\left(t\right)=\left[u_{A 1}\left(t\right), u_{A 2}\left(t\right) \right]=\sum\limits_{i=0}^{p}\overrightarrow{u}_{A}^{(i)}\left(t\right)$ where unknown vector-functions $\overrightarrow{u}_{A}^{(i)}\left(t\right)$ are the solutions of the following Cauchy problems
\begin{equation}\label{Chapt_3_5_eq_13}
\begin{array}{c}
  \cfrac{d}{dt}\overrightarrow{u}^{(0)}_{A}\left(t\right)=-E^{(2)}\overrightarrow{u}^{(0)}_{A}\left(t\right)+\overrightarrow{\phi}\left(\xi\right), \quad \overrightarrow{u}^{(0)}_{A}\left(0\right)=\left[0, 1\right]^{T},\\
  \cfrac{d}{dt}\overrightarrow{u}^{(i)}_{A}\left(t\right)=-E^{(2)}\overrightarrow{u}^{(i)}_{A}\left(t\right)+\cfrac{1}{i!}{\bigg.\cfrac{d^{i}}{d \tau^{i}}\left(\tau \mathfrak{N}_{1}\left(\sum\limits_{i=0}^{\infty}\tau^{i}\overrightarrow{u}_{A}^{(i)}\left(\xi\right)\right)\right)\bigg|_{\tau=0}}=\\
\end{array}
\end{equation}
$$\begin{array}{c}
  =-E^{(2)}\overrightarrow{u}^{(i)}_{A}\left(t\right)+A_{i-1}\left(\mathfrak{N}_{1}\left(\cdot\right), \left[\overrightarrow{u}^{(k)}_{A}\left(\xi\right)\right]_{k=0}^{i-1}\right), \quad \overrightarrow{u}^{(i)}_{A}\left(0\right)=0,\quad i\in\N. \\
\end{array}
$$
The results shown on figure \ref{image2} tell us that the  ADM \eqref{Chapt_3_5_eq_13} is divergent on the segment $\left[0, 2\right].$

Applying the FD-method to the Cauchy problem \eqref{Chapt_3_5_eq_1}, \eqref{Chapt_3_5_eq_2} we used the uniform grid with the step $h=1/5.$ We have found five approximations $\overset{p}{\overrightarrow{u}}\left(t\right),$ $p\in\overline{0, 4}$ of the exact solution on the segment $\left[0, 6\right].$ Also we have calculated the absolute error values $$\delta_{i}^{(p)}\left(t\right)=\left|\overset{p}{u}_{i}\left(t\right)-u_{i}\left(t\right)\right|,\quad i=1,2,\quad p\in \overline{0, 4},$$
$$\delta^{(p)}\left(t\right)=\left\|\overset{p}{\overrightarrow{u}}\left(t\right)-\overrightarrow{u}\left(t\right)\right\|,\quad p\in \overline{0, 4}.$$

 The numerical results showed on figure \ref{image3} confirm that the FD-method for the Cauchy problem \eqref{Chapt_3_5_eq_1}, \eqref{Chapt_3_5_eq_2} converges with the exponential rate to the exact solution of the problem.

\pagebreak

\end{document}